\documentclass[11pt]{article}

\ifx\comments\undefined
  \usepackage[letterpaper,margin=1in,includefoot]{geometry}
  \newcommand{\XXXcomment}[1]{}
  \newcommand{\XXXcommentR}[1]{}
  \newcommand{\XXXcommentG}[1]{}
\else
  \usepackage[letterpaper,margin=.5in,right=2.5in,marginpar=2in,includefoot]{geometry}
  \newcommand{\XXXcomment}[1]{\marginpar{\color{blue}{\footnotesize #1}}}
  \newcommand{\XXXcommentR}[1]{\marginpar{\color{red}{\footnotesize #1}}}
  \newcommand{\XXXcommentG}[1]{\marginpar{\color{green}{\footnotesize #1}}}
\fi

\ifx\arXiv\undefined
  \usepackage[pdftex,
              bookmarks=true,
              bookmarksnumbered=true,
              pdfborder={0 0 0},
              plainpages=false,
              unicode,
              pdfpagelabels]{hyperref}
\fi

\usepackage{amssymb,amsthm,latexsym,amsmath,mathrsfs,stmaryrd,graphicx,fancyhdr,paralist,longtable, multirow, wasysym}

\usepackage[x11names]{xcolor}
\usepackage{colortbl}

\usepackage[compact]{titlesec}
\usepackage[titletoc,title]{appendix}
\usepackage{multicol}

\usepackage{bussproofs}

\usepackage{linguex}

\ifx\arXiv\undefined
  \usepackage[sort]{natbib}
\fi

\usepackage{caption,float}
\usepackage{tikz}
\usetikzlibrary{trees,arrows,automata,shapes,decorations,decorations.pathmorphing}

\usepackage[normalem]{ulem}

\usepackage[makeroom]{cancel}

\theoremstyle{definition} 
\newtheorem{definition}{Definition}[section]
\newtheorem{theorem}[definition]{Theorem}

\newtheorem{corollary}[definition]{Corollary}
\newtheorem{example}[definition]{Example}

\newtheorem{notation}[definition]{Notation}
\newtheorem{fact}[definition]{Fact}

\definecolor{darkgray}{rgb}{.1, .1, .14}
\newenvironment{myproof}{\footnotesize\begin{proof}\color{darkgray}}{\end{proof}\normalsize}
\newenvironment{myproofbis}{\footnotesize\begin{proof}[The proof of the primed version\nopunct]\color{darkgray}}{\end{proof}\normalsize}
\newtheorem{thm}{Theorem}
\newenvironment{thmbis}[1]
  {\renewcommand{\thethm}{\ref{#1}$'$}%
   \addtocounter{thm}{-1}%
   \begin{thm}}
  {\end{thm}}

\newcommand{\ourtitle}{Suszko's Problem: Mixed Consequence and Compositionality}

\newcommand{\paul}{Paul Egr\'e}

\newcommand{\paulinstitute}{IJN / ENS}

\newcommand{\emmanuel}{Emmanuel Chemla}

\newcommand{\emmanuelinstitute}{LSCP / ENS}

\newcommand{\emmanuelFunding}{Funded by the European Research Council under the European Union's Seventh Framework Programme (FP/2007-2013) / ERC Grant Agreement n. 313610 and by grant ANR-10-IDEX-0001-02.}
  
\newcommand{\paulFunding}{Funded by the ANR project ``Trivalence and
  Natural Language Meaning'' (ANR-14-CE30-0010), by grant ANR-10-IDEX-0001-02, and by Ministerio de Econom\'ia y Competitividad, Gobierno de Espana as part of the project "Logic and Substructurality", grant number FFI2017-84805-P.}

\ifx\arXiv\undefined
\hypersetup{ 
  pdfauthor={\emmanuel{}, \paul{}}, 
  pdftitle={\ourtitle{}}
}
\fi

\author{
  \emmanuel{}\footnote{\emmanuelFunding{}}\\
  \small\emmanuelinstitute{}
  \and
  \paul{}\footnote{\paulFunding{}}\\
  \small\paulinstitute{}
 }

\date{}

\def\truthrelation{\mathrel|\joinrel\equiv}

\usepackage[T1]{fontenc}
\usepackage[sc,osf]{mathpazo}

\usepackage{microtype}

\newcommand{\interp}[1]{\llbracket #1 \rrbracket}

\newcommand{\dint}[1]{{\left\vert\kern-0.25ex\left\vert\kern-0.25ex\left\vert #1 
    \right\vert\kern-0.25ex\right\vert\kern-0.25ex\right\vert}}
    \newcommand{\sem}[1]{\mbox{$[\![\![ #1 ]\!]\!]$}}

\begin{document}

\renewcommand*{\thefootnote}{\fnsymbol{footnote}}
\noindent {\Large 
Suszko's problem: mixed consequence and {compositionality}%
\footnote{{Acknowledgements: We are grateful to Johan van Benthem, Denis Bonnay, Pablo Cobreros, Andreas Fjellstad, Ole Hjortland, Jo{\~a}o Marcos, Francesco Paoli, David Ripley, Robert van Rooij, Lorenzo Rossi, Philippe Schlenker, Benjamin Spector, Shane Steinert-Threlkeld, Heinrich Wansing for various exchanges and helpful suggestions. We are particularly grateful to David Ripley and Jo{\~a}o Marcos for their encouragements and for enlightening discussions regarding each of their respective works on this topic. We thank Andrew Arana and two anonymous referees for valuable comments that helped improve the paper. We also thank audiences in workshops held in Pamplona, Munich, Paris, and Amsterdam. This research was partially supported by the program "Logic and Substructurality", grant number FFI2017-84805-P, from the Ministerio de Econom\'ia, Industria y Competitividad, Government of Spain, by the French ANR program TRILOGMEAN (ANR-14-CE30-0010), by the European Research Council under the European Union's Seventh Framework Programme (FP/2007-2013) / ERC Grant Agreement n. 313610, and by grant ANR-10-IDEX-0001-02 for research carried out at the DEC-ENS in Paris.}
}}\\
\renewcommand*{\thefootnote}{\arabic{footnote}}
\setcounter{footnote}{0}
Emmanuel Chemla$^\textrm{a}$ \& Paul Egr\'e$^\textrm{b}$ \\
{\scriptsize a. LSCP, D\'epartement d'\'etudes cognitives, ENS, EHESS, CNRS, PSL University, 75005 Paris, France}\\
{\scriptsize b. Institut Jean Nicod, D\'epartement d'\'etudes cognitives \& D\'epartement de philosophie, ENS, EHESS, CNRS, PSL University, 75005 Paris, France}

\begin{abstract}

\noindent Suszko's problem is the problem of finding the minimal number of truth values needed to semantically characterize a syntactic consequence relation. Suszko proved that every Tarskian consequence relation can be characterized using only two truth values. Malinowski showed that this number can equal three if some of Tarski's structural constraints are relaxed. By so doing, Malinowski introduced a case of so-called \emph{mixed consequence}, allowing the notion of a designated value to vary between the premises and the conclusions of an argument. In this paper we give a more systematic perspective on Suszko's problem and on mixed consequence. First, we prove general representation theorems relating structural properties of a consequence relation to their semantic interpretation, uncovering the semantic counterpart of substitution-invariance, and establishing that (intersective) mixed consequence is fundamentally the semantic counterpart of the structural property of monotonicity. We use those theorems to derive maximum-rank results proved recently in a different setting by French and Ripley, as well as by Blasio, Marcos and Wansing, for logics with various structural properties (reflexivity, transitivity, none, or both). We strengthen these results into exact rank results for \emph{non-permeable} logics (roughly, those which distinguish the role of premises and conclusions). We discuss the underlying notion of rank, and the associated reduction proposed independently by Scott and Suszko. As emphasized by Suszko, that reduction fails to preserve compositionality in general, meaning that the resulting semantics is no longer truth-functional. We propose a modification of that notion of reduction, allowing us to prove that over compact logics with what we call \emph{regular connectives}, rank results are maintained even if we request the preservation of {truth-functionality} and additional semantic properties.

\end{abstract}

\bigskip

\section{Suszko's problem}

This paper deals with an issue famously raised by the Polish logician Roman Suszko concerning many-valuedness in logic. The issue, which we call \emph{Suszko's problem}, may be presented as follows: given a language and a consequence relation between the formulae of that language, satisfying some general structural constraints, what is the minimum number of truth values required to characterize the relation semantically? 

Suszko gave a precise answer to that problem for a large class of consequence relations. He proved that every Tarskian relation, namely every reflexive, monotonic, and transitive consequence relation, can be characterized by a class of two-valued models \citep{suszko1977fregean, tsuji1998many, wansing2008suszko, caleiro2012many, caleiro2015bivalent}. Suszko's problem and his result are of fundamental importance to understand the role and nature of truth values in logic. What Suszko had in mind is that even as a language is semantically interpreted over more than two truth values, only two may suffice to semantically characterize logical consequence. Based on Suszko's result, it has become customary to distinguish two roles for truth values: a referential role, pertaining to the semantic status of specific formulae, and an inferential role, concerned with the logical relations between them (see \citealp{shramko2011truth} for an extensive review). The referential role is encoded by what we may call \emph{algebraic truth values}, those introduced by the semantics, whereas the inferential role is encoded by what we may call \emph{logical truth values}, which may be thought of as more abstract entities \citep{suszko1975remarks}. Underpinning Suszko's theorem is indeed the distinction between designated and undesignated algebraic values: Suszko's theorem can be understood in relation to the semantic definition of logical consequence, whereby a conclusion follows from a set of premises provided it is not the case that the conclusion is undesignated when the premises are all designated. 

Suszko extracted a bolder lesson from his result: he surmised that logical consequence is essentially two-valued. This claim, known as \emph{Suszko's thesis}, depends for its correctness on an assessment of what to count as a consequence relation, and in particular on whether every {bona fide} consequence relation should obey the Tarskian constraints. In response to Suszko's thesis, Malinowski showed that if some of Tarski's constraints are relaxed, then logical consequence may in some cases be essentially three-valued. To show this, \cite{malinowski:q} introduced the notion of $q$-consequence. Semantically, $q$-consequence can be described as a variety of \emph{mixed consequence}, allowing the notion of a designated value to vary between the premises and the conclusion of an argument. In three-valued terms, assuming the values $1$, $0$ and $\frac{1}{2}$, $q$-consequence may be defined as follows: whenever the premises of an argument do not take the value $0$, the conclusion must take the value $1$. What Malinowski's result shows is that this notion of consequence (aka $ts$-consequence, see \citealp{tcs}) is not reducible to fewer than three values. $ts$-consequence is monotonic and transitive, but it is not reflexive. Other notions of mixed consequence have since been discussed in the literature, in particular Frankowski's $p$-consequence (\citealp{frankowski2004}), whose semantic counterpart is $st$-consequence, symmetric to $ts$-consequence, requiring the conclusion to not take the value $0$ when the premises all take the value $1$. Unlike $ts$-consequence, $st$ is reflexive, but it is not transitive. Like $ts$-consequence, $st$-consequence is logically three-valued.

Recently, \cite{french2017valuations}, {and \cite{blasio2017inferentially}}, have derived a variety of Suszko-type results, providing a very systematic perspective on Suszko's problem.{\footnote{{Our work was developed independently of the work of \cite{blasio2017inferentially}, and only referred to the work of \cite{french2017valuations} at the time we submitted the first version of this paper. We are grateful to J. Marcos for bringing his joint paper with C. Blasio and H. Wansing to our notice. We do not undertake a detailed comparison of the three frameworks here: suffice it to point out that they differ both in style and in what they assume as the basic notions.}}} They prove that every monotonic logic is at most 4-valued, that every monotonic reflexive logic is at most 3-valued, likewise that every monotonic transitive logic is at most 3-valued, and finally that every logic satisfying all of Tarski's constraints is at most 2-valued. In this paper, we pursue a very similar goal, and in particular we derive the same results, but within a different framework (much in line with \citealp{sep-logic-algebraic-propositional, bloom1970theorems, andreka2001algebraic, font2003generalized, wojcicki1988theory}), and coming from a different angle (see \citealp{CES2016}). Like French and Ripley, {and like Blasio, Marcos and Wansing}, we will take advantage of systematic correspondences between the structure of arguments and their semantic interpretations. {Our perspective and contribution differ in three ways.}

First, the notion of mixed consequence is explicitly given a central role in our inquiry, both as a topic of study and as a technical tool.
In section \ref{sec:correspondence}, in particular, we prove that intersections of mixed consequence relations constitute an exact counterpart of the structural property of monotonicity. This fact, although implicit in the precursor works of \cite{humberstone1988heterogeneous} and indeed of \cite{french2017valuations}, {and \cite{blasio2017inferentially}}, was not brought to the fore in exactly that form (see in particular \citealp{CES2016} for a closely related result, but missing the simpler characterization). 
Secondly, we strengthen these results into \emph{exact} rank results for \emph{non-permeable} logics, roughly, those which distinguish the role of premises and conclusions: reflexive and transitive logics are of rank exactly 2, logics which are either reflexive or transitive but not both are of rank exactly 3, all others are of rank 4.
Third, we give a specific attention to the issue of semantic compositionality. \cite{suszko1977fregean} pointed out that his reduction technique generally ends up in a semantics {that fails compositionality, understood as truth-functionality (see \citealt{beziau2001sequents, font2009taking, marcos2009non, caleiro2012many, caleiro2015bivalent} for discussions, and our refinements of the notion of compositionality below)}. This failure, as stressed by \cite{font2009taking}, makes it ``difficult to call it a semantics, which makes the meaningfulness of his Reduction hard to accept''. {One of our goals in this paper is to clarify this limitation and go beyond it. {W}e will extend the notion of reduction to prove that every compact logic made of connectives satisfying special properties {(\textit{regular connectives})} can be given a {truth-functional} semantics of the same 2/3/4 ranks depending on their structural properties.}

Our work is organized as follows. In section 2, we introduce our framework and natural structural constraints usually imposed on a logic (substitution-invariance, monotonicity, reflexivity, and transitivity), as well as natural semantic constraints on the interpretation of a logic (e.g., compositionality, truth-functionality).
In section 3, we prove representation theorems establishing the tight mapping between structural properties and semantic properties. In the first part of section 4, we then apply the Scott-Suszko reduction technique to derive analogues of the Suszko results obtained by French and Ripley, {as well as Blasio, Marcos and Wansing}. 
In the rest of section 4, we turn to the problem of {truth-functionality}. We introduce a novel type of reduction, namely a construction that transforms a semantics into another semantics with fewer truth values, while preserving both the consequence relation and {truth-functionality}. We call this the \emph{truth-{functional}} Scott-Suszko reduction, allowing us to isolate a natural class of logics {-- compact regular logics--} over which the resulting semantics is {truth-functional}. This modification is substantive: it brings further {significance} to Suszko's thesis by proving not only that the notion of logical many-valuedness can be meaningful beyond two values, but also that reduction of many values to a minimum is not necessarily a threat on {truth-functionality}. 
In effect, we prove that under appropriate circumstances, the truth values needed to preserve the consequence relation are sufficient to preserve truth-{functionality} as well. 

\section{Framework}\label{sec:framework}

This section lays out the basic ingredients for the rest of this paper. We define a logic in syntactic terms, and right away highlight the main structural constraints on consequence relations isolated by Tarski and the Polish school (see \citealp{tarski1930, bloom1970theorems, wojcicki1973matrix, wojcicki1988theory}), namely: substitution-invariance, monotonicity, reflexivity and transitivity. Those properties play a central role in Suszko's reduction. We then articulate what we mean by a semantics for a logic, and state various semantic properties which we will prove to be counterparts of those structural properties in section \ref{sec:correspondence}.

\subsection{Syntactic notions}

\subsubsection{Logic}
	
{Throughout the paper, we work within a multiple-premise multiple-conclusion logic, following the tradition of \cite{gentzen1935investigations}, \cite{scott1974completeness} and \cite{shoesmith1978multiple}. One of the reasons for doing so is simplicity in that premises and conclusions are handled symmetrically. This choice is not neutral and some of our results are likely to differ in a multi-premise single-conclusion setting. We leave that issue aside, and define a logic as follows:}

\begin{definition}[Logic]
A \emph{logic} is a triple $\langle\mathcal{L}, \mathcal{C}, \vdash\rangle$, with 
$\mathcal{L}$ a language (set of formulae), $\mathcal{C}$ a (possibly empty) set of connectives,
and $\vdash$ a \emph{consequence relation}, what we call a \emph{formula-relation} (i.e. a subset of $\mathcal{P}(\mathcal{L})\times\mathcal{P}(\mathcal{L})$, also called a set of \emph{arguments}, where each argument is a relation between sets of formulae). 
\end{definition}

\begin{definition}[Sentential Language]

A language is called \emph{sentential} if it is obtained in the usual way from a set $\mathcal{A}$ of atomic formulae (a subset of $\mathcal{L}$) and a set $\mathcal{C}$ of \emph{connectives}, that is formula-functions of type $c:\mathcal{L}^n\longrightarrow\mathcal{L}$. A logic is called sentential if it is defined over a sentential language.
\end{definition}

\noindent We will only consider sentential languages and sentential logics in what follows.

\begin{definition}[Substitution] A substitution $\sigma$ for a sentential language is an endomorphism of $\mathcal{L}$, such that for every formula $F(p_{1},...,p_{n})$ constructed from the atoms $p_i$, $F[\sigma]=F(\sigma(p_{1}),...,\sigma(p_{n}))$. 
\end{definition}

\subsubsection{Fundamental structural properties}

\begin{definition}[Substitution-invariance]
A consequence relation $\vdash$ is \emph{substitution-invariant} iff for every substitution $\sigma$:
$$\textrm{If } \Gamma \vdash \Delta \textrm{, then } \Gamma[\sigma] \vdash \Delta[\sigma]$$
\end{definition}

\begin{definition}[Monotonicity]
A consequence relation $\vdash$ is \emph{monotonic} if:
$$\forall \Gamma_1\subseteq\Gamma_2, \Delta_1\subseteq\Delta_2:
\Gamma_1 \vdash \Delta_1
\textrm{ implies } 
\Gamma_2 \vdash \Delta_2.$$

\end{definition}

	\begin{definition}[Reflexivity]
A consequence relation $\vdash$ is \emph{reflexive} if 
for every formula $F$:
$$ F \vdash F$$
\end{definition}

\begin{definition}[Transitivity]

	A consequence relation $\vdash$ is \emph{transitive} iff: 
	
	$$\textrm{if }\Gamma\not\vdash \Delta\textrm{, then there are }
	\Gamma' \supseteq\Gamma, \Delta'\supseteq\Delta\textrm{ such that }\Gamma'\not\vdash \Delta'\textrm{ and }\Gamma'\cup\Delta'=\mathcal{L}.$$
	\end{definition}

\noindent There are numerous ways of defining transitivity, in particular in multi-conclusion settings. They are extensively reviewed in \cite{ripley2017transitivity}. The definition proposed above is the strongest definition from this survey, this is the one that will prove relevant in the next section.%
\footnote{Let us show why the usual notion of transitivity follows.
Consider that $A\not\vdash C$. The above notion entails that there are $\Gamma,\Delta$ such that $B$ is in one of them and $A, \Gamma \not\vdash C, \Delta$. It follows by monotonicity that either $B \not\vdash C$ ($B\in\Gamma$) or that $A \not\vdash B$ ($B\in\Delta$). By contraposition we obtain that $A\vdash B$ and $B \vdash C$ entail $A\vdash C$.}

\begin{definition}
We say that a logic is \emph{substitution-invariant} / \emph{monotonic} / \emph{reflexive} / \emph{transitive} if it has a substitution-invariant / monotonic / reflexive / transitive consequence relation, respectively.
\end{definition}

\subsubsection{A further structural constraint: non-permeability}

We have introduced structural properties from the Tarskian tradition, which are supposed to delineate what an appropriate consequence relation may be. These properties, however, do not exclude pathological consequence relations, such as the null and the universal consequence relations, or relations for which there would be no non-empty $\Gamma$ and $\Delta$ such that $\Gamma\not\vdash\Delta$. 
We introduce here a formal property that a consequence relation should \emph{not} have:

\begin{definition}[Permeability]\label{def:permeable}
A consequence relation is \emph{permeable} if it is \emph{left-to-right} or \emph{right-to-left permeable}, in the following sense:
\begin{description}
\item[Left-to-right permeability:]
	$\forall\Gamma,\Delta,\Sigma: \Gamma,\Sigma\vdash\Delta \Rightarrow \Gamma\vdash \Sigma,\Delta$
\item[Right-to-left permeability:]
	$\forall\Gamma,\Delta,\Sigma: \Gamma\vdash\Sigma,\Delta \Rightarrow \Gamma, \Sigma\vdash\Delta$
\end{description}
By extension, a logic is called \emph{permeable} if its consequence relation is permeable.
\end{definition}

If a logic is not permeable, then its consequence relation is neither universal nor trivial (when negated, either constraint implies non-triviality as per the antecedent and non-universality as per the denied consequent). We may describe a permeable consequence relation as one that would confuse the role of premises with the role of conclusions, or the reverse. But we need to introduce more semantic machinery before we describe the full virtue of this constraint (see Theorem~\ref{thm:permeablepolarized}).

\subsection{Semantic notions}

\subsubsection{Semantics}

We call a \emph{semantics} a structure that interprets each of the components of a logic, namely formulae, connectives, and the consequence relation. The natural semantic counterpart for formulae are \emph{propositions}, understood as functions from worlds to truth values.
The semantic counterpart of a connective is therefore a function from propositions to propositions, and the semantic counterpart of a consequence relation is a relation between sets of propositions. But functions and relations over sets of propositions may be obtained from functions and relations over sets of truth values: just like connectives may be interpreted through \emph{truth}-functions, we show that consequence relations may be interpreted algebraically through what we call \emph{truth-relations}.
An overview of the definitions to be introduced is given in Table~\ref{table:terminology}.

\subsubsection*{Semantics with propositions}

\begin{definition}[Semantics]
A \emph{semantics} for a logic is a triple $\langle\mathcal{V}, \mathcal{W}, \interp{\_}\rangle$, with 
\begin{itemize}
\item $\mathcal{V}$ a set of \emph{truth values}, 
\item $\mathcal{W}$ a set of \emph{worlds}, 
\item $\interp{\_}$ an \emph{interpretation function}, which associates 
 	to every formula $F$ in $\mathcal{L}$ a \emph{proposition} $\interp{F}$, defined as a function from worlds to truth values, 
	to every $n$-ary connective $c$ in $\mathcal{C}$ a function $\interp{c}$ from $n$-tuples of propositions to propositions, and 
	to $\vdash$ it associates a relation $\interp{\vdash}$ between sets of propositions (henceforth $\models$).
\end{itemize}
\end{definition}

Semantics can match syntax to a variable extent. Below are two main properties which pertain both to the connectives and to the consequence relation:

\begin{definition}[Compositionality, Soundness-and-completeness, Adequacy]\label{def:semprop} Given a logic and a semantics:
\begin{itemize}
\item {The interpretation of a}n $n$-ary connective $c$ is \emph{compositional} if for all formulae $F_{1},...,F_{n}$
	$$\interp{c(F_1, ..., F_n)} = \interp{c}(\interp{F_1}, ..., \interp{F_n})$$
\item {The interpretation of t}he consequence relation $\vdash$ is \emph{sound and complete} if
for all sets of formulae $\Gamma, \Delta$
	$$\Gamma \vdash \Delta
	\textrm{ iff }
	\interp{\Gamma}\models \interp{\Delta}$$ 
	
\item A semantics for a logic is called 
\emph{compositional} if it interprets every connective in a compositional way, 
\emph{sound and complete} if it interprets the consequence relation in a sound and complete way, and 
\emph{adequate} if it has both properties. 
\end{itemize}
\end{definition}

\subsubsection*{Semantics from truth values}

So far, the consequence relation and the connectives have been interpreted at the level of propositions. But a semantics may intrinsically follow from a definition of connectives and consequence relations at the level of truth values, and uniformly so world by world:

\begin{definition}[Truth-Functionality, Truth-Relationality, Truth-Interpretability]\label{def:semtruth} Given a logic and a semantics:
\begin{itemize}

\item {The interpretation of a}n $n$-ary connective $c$ is \emph{truth-functional} if 
	there is a \emph{truth-function} $\sem{c}$, 
	that is a function from $n$ truth values to truth values (of type $\mathcal{V}^n\to\mathcal{V}$) such that:
	$$
	\textrm{for all} \textrm{ propositions }P_1, ..., P_n: \forall v\in\mathcal{W}: 
		\interp{c}(P_1, ..., P_n)(v)=\sem{c}(P_1(v), ..., P_n(v))
	$$

\item {The interpretation of t}he consequence relation $\vdash$ is \emph{truth-relational} if 
	there is a \emph{truth-relation} $\sem{\vdash}$, also noted $\truthrelation$, 
	that is a relation between sets of truth values (of type $\mathcal{P}(\mathcal{V})\times\mathcal{P}(\mathcal{V})$) such that:
	$$
	\textrm{for all sets of propositions } \mathcal{S}, \mathcal{P}:
			\mathcal{S}\models\mathcal{P} 
		\textrm{ iff }
		\forall v\in\mathcal{W}: \mathcal{S}(v)\truthrelation\mathcal{P}(v)
	$$

\item 
A semantics for a logic is called \emph{truth-functional} if it interprets every connective in a truth-functional way,
\emph{truth-relational} if it interprets the relation $\vdash$ in a truth-relational way,
\emph{truth-interpretable} if it has both properties.
\end{itemize}
\end{definition}

\begin{theorem} To validate the previous definition, and in particular the notations $\sem{C}$ and $\sem{\vdash}$/$\truthrelation$, we state that (for $\mathcal{W}\not=\emptyset$, which is what we are interested in):
\begin{itemize}
\item The truth-function $\sem{C}$ associated to a truth-functional connective $C$ is uniquely determined by $\interp{C}$.
\item The truth-relation $\truthrelation$ associated to a truth-relational consequence relation $\vdash$ is uniquely determined by $\models$.
\end{itemize}
\end{theorem}

\begin{myproof}
Truth-functionality and truth-relationality alike are about the fact that $\sem{\_}$ is uniform across worlds. The relevant truth-functions and truth-relation are simply that constant function (value) or relation taken in each of the worlds. To put it differently, a truth-functional connective associates a constant proposition as the output of its application to constant propositions; the truth-function ought to preserve this and therefore is determined by:
$\sem{c}(\alpha_1, ..., \alpha_n)$ is the value taken by the constant proposition $\interp{c}(\overline{\alpha_1}, ..., \overline{\alpha_n})$, where $\alpha_1, ..., \alpha_n$ are truth values and $\overline{\alpha_1}, ..., \overline{\alpha_n}$ constant propositions taking these values. Similarly, the truth-relation of a truth-relational consequence relation is determined by $\gamma\truthrelation\delta$ iff $\overline{\gamma} \models \overline{\delta}$, with $\gamma, \delta$ sets of truth values and $\overline{\gamma}, \overline{\delta}$ the corresponding sets of constant propositions.
\end{myproof}

Table~\ref{table:terminology} summarizes the semantic definitions introduced so far, in parallel for connectives and relations, first at the level of propositions (Definition~\ref{def:semprop}) and then at the level of truth values (Definition~\ref{def:semprop}). In this Table we also introduce names for conjunctions of properties: 

\begin{definition}
A compositional and truth-functional semantics is called \emph{truth-compositional}; a sound and complete and truth-relational semantics is called \emph{truth-(sound and complete)}; an adequate and truth-interpretable semantics is called \emph{truth-adequate}.
\end{definition}

\begin{table}[!h]

\begin{centering}\framebox{\begin{tabular}{l|ll|l}
& \emph{Proposition level} & \emph{Truth-value level} & \emph{Conjunctively}\\ \hline
\emph{Connectives} & compositional & truth-functional & truth-compositional \\
\emph{Consequence} & sound and complete & truth-relational & truth-(sound and complete)\\ \hline
\emph{Conjunctively} & adequate & truth-interpretable & truth-adequate\\
\end{tabular}}
\caption{Summary of the properties introduced in parallel for connectives and consequence relations (or both, when combined at the level of their logic). These properties belong to the level of propositions, or truth values (or both, if we consider the conjunction of the properties at the proposition level and its associate at the truth-value level).}
\label{table:terminology}
\end{centering}
\end{table}

\noindent {Note that according to our definitions, it would be possible} for a semantics to be truth-functional without being compositional {proper (for example, it could happen that for all $v$, $(\interp{\neg}\interp{F})(v)=\sem{\neg}(\interp{F}(v))$, but $\interp{\neg}\interp{F}\neq \interp{\neg F}${, simply note that the first truth-functional equality does not involve $\neg F$})}. When the semantics is defined directly at the truth-value level, however, truth-functionality directly implies compositionality (and therefore coincides with truth-compositionality), {as the following shows}.

\begin{fact}\label{fact:semanticsfromtruth}
An efficient, and usual way to determine a truth-functional semantics is to define the interpretation function for formulae only on atoms and extend it to all formulae by defining directly the truth-functions $\sem{C}$ for all connectives: 
$$\interp{C(F_1, ..., F_n)}\textrm{ is defined as }\forall v: \interp{C(F_1, ..., F_n)}(v):=\sem{C}(\interp{F_1}(v), ..., \interp{F_n}(v))$$

Similarly, a truth-relational semantics may be obtained 
by choosing a truth-relation~$\truthrelation$ and then deriving the interpretation of the consequence relation $\models$ as: 
$$\interp{\Gamma}\models\interp{\Delta} \quad\textrm{ iff (def) }\quad \forall v: \interp{\Gamma}(v)\truthrelation\interp{\Delta}(v)$$
\end{fact}

\subsubsection*{Truth-relationality}

We pause on the notion of truth-relationality, which is the least standard (see \citealp{CES2016}). First, we note that in the absence of other constraints, the truth-relationality of \emph{some} sound and complete semantics for a logic is always guaranteed. In fact, the conjunction of all properties above can always be achieved:

\begin{fact}\label{fact:t-adeqforall}
Every {sentential} logic has a truth-adequate semantics.
\end{fact}
		
\begin{myproof}
Given a logic, 
	let 
	$\mathcal{V}:=\mathcal{L}$, 
	$\mathcal{W}:=\{v_{0}\}$, 
	$\interp{F}({v}_{0}):=F$,
	and let $\sem{c}$ be $c$ itself for every connective and $\truthrelation$ be $\vdash$ itself. 
By construction, $\langle \mathcal{V}, \mathcal{W}, \interp{\_} \rangle$ is a truth-interpretable semantics, 
	compositional ($\interp{C}(\interp{F_1}, ..., \interp{F_n})=C(F_1, ..., F_n)=\interp{C(F_1, ..., F_n)}$)
	and sound and complete ($\Gamma \vdash \Delta$ iff for all $v\in \mathcal{W}$, $\interp{\Gamma}(v) \truthrelation \interp{\Delta}(v)$).
\end{myproof}

Despite that result, truth-relationality is not unsubstantial, and it certainly does not follow from the other properties a semantics may have.

\begin{fact}\label{fact:supervaluationism}
Not every adequate semantics that is truth-functional is truth-relational.
\end{fact}

\begin{myproof}

Many degenerated examples could be used to demonstrate the result, but let us consider the richer case of so-called \emph{supervaluationism} (viz. \citealp{cobreros:mixed}). Consider the language $\mathcal{L}$ of standard propositional logic over a denumerable set of atoms $\mathcal{A}$ with set of connectives $\mathcal{C}=\{\neg, \vee, \wedge\}$. We can define the so-called supervaluationist semantics on $\mathcal{L}$, which is exactly like the classical bivalent strong semantics over $\mathcal{L}$, {with $\mathcal{W}$ as the set of all valuations (see \ref{def:val})}, except that the relation is interpreted as:
$$\Gamma\models_{\textrm{\textsc{sv}}}\Delta 
	\quad\textrm{ iff (def) }\quad
	\forall\mathcal{I}\subseteq\mathcal{W}:
	(\forall F_p\in\Gamma: \forall v\in\mathcal{I}: \interp{F_p}(v)=1)
	\Rightarrow 
	(\exists F_c\in\Delta: \forall v\in\mathcal{I}: \interp{F_c}(v)=1)$$

Supervaluationist consequence induces a supervaluationist $\vdash_{\textrm{\textsc{sv}}}$ syntactic consequence relation over $\mathcal{L}$ ($\Gamma \vdash_{\textrm{\textsc{sv}}} \Delta$ iff (def) $\interp{\Gamma}\models_{\textrm{\textsc{sv}}}\interp{\Delta}$). The supervaluationist semantics is sound and complete by construction, it is compositional and truth-functional just like classical logic. But there could be no truth-relation. The reason is as follows:
\begin{itemize}
\item $\not \models_{\textrm{\textsc{sv}}} \interp{p}, \interp{\neg p}$. To prove this, call the (empty) premise set $\Gamma$ and the conclusion set $\Delta$. Even though
	$\forall F_p\in\Gamma: \forall v\in\mathcal{W}: \interp{F_p}(v)=1$ (trivially satisfied because $\Gamma=\emptyset$),
	it is not the case that $\exists F_c\in\Delta: \forall v\in\mathcal{W}: \interp{F_c}(v)=1$.

\item $\models_{\textrm{\textsc{sv}}} \interp{p\vee\neg p}, \interp{p\wedge\neg p}$, this is because $\interp{p\vee\neg p}$ is the constant proposition with value $1$.

\item For all valuations, $\{\interp{p}(v), \interp{\neg p}(v)\} = \{0,1\} = \{\interp{p\wedge \neg p}(v), \interp{p\vee \neg p}(v)\}$.
\end{itemize}
The three facts above are incompatible with the possibility of a truth-relation, which would sometimes have to hold for $\{0,1\}$ in its second argument and sometimes not.
\end{myproof}

\noindent {The} semantics {we proposed for supervaluationism} is therefore compositional and truth-functional, but not truth-relational. The logic it corresponds to may have a truth-relational semantics, but this may require another set of truth values {and another set of worlds}. 
A similar example can be found in \cite{CES2016}, with a discussion of a logic and a natural semantics for it that is adequate, truth-functional but not truth-relational, namely the trivalent semantics $ss\cup tt$ for, e.g., some propositional language (the union of strong Kleene consequence and LP consequence).

\subsubsection{Fundamental semantic properties}

We now state several fundamental properties of a semantics. 
We will refer to these properties as semantic properties for short. 
In the next section, we will prove that they are parallel to the structural constraints we stated for a logic (an overview of the parallelism is given in Table~\ref{tab:representation}).
\begin{definition}\label{def:val}
A semantics for a {sentential} logic is called \emph{valuational} if its set of worlds is the set of \emph{valuations}, functions from atoms to truth values, and if for all atoms $p$ and all valuations $v$, $\interp{p}(v)=v(p)$.
\end{definition}

\begin{fact}\label{fact:truth-compositionalityandvaluationality}
For a truth-compositional semantics
it is also appropriate (and will be useful) to interpret worlds as valuations: the truth value assigned to a formula in a given world is entirely determined by the truth values assigned to the atoms (involved in the formula) in that world, by way of the successive application of the truth-functions of the connectives that make up the formula. As a consequence, one may drop redundant worlds/valuations that would assign the same values to all atoms.
Truth-compositionality does not entail valuationality, however, because $\mathcal{W}$ may not cover the whole set of possible valuations.
\end{fact}

\begin{notation}\label{not:val}
For a valuational semantics $\langle\mathcal{V}, \mathcal{V}^{\mathcal{A}}, \interp{\_}\rangle$, we may drop from the notation the set of worlds $\mathcal{W}=\mathcal{V}^{\mathcal{A}}$ and refer to it as $\langle\mathcal{V}, \interp{\_}\rangle$. Furthermore, we may refer to $\interp{F}(v)$ as $v(F)$, by extension of the notation for atomic formulae, whenever the interpretation of worlds as valuations is warranted (for valuational semantics or truth-compositional semantics as in Fact~\ref{fact:truth-compositionalityandvaluationality}).
\end{notation}

\begin{definition}[Strong semantics]
We say that a semantics is \emph{strong} if it is truth-adequate and valuational (all the properties above).
\end{definition}

\begin{definition}[Mixed truth relation]\label{def:mixed}
A {truth-relation} $\truthrelation$ over a set of truth values $\mathcal{V}$ is called \emph{mixed} if there are two sets $\mathcal{D}_p$ and $\mathcal{D}_c$ included in $\mathcal{V}$ such that for all sets of values $\gamma$ and $\delta$:
 $$\gamma \truthrelation \delta\ \textrm{iff}\ \gamma \subseteq \mathcal{D}_p\ \textrm{entails}\ \delta \cap \mathcal{D}_c \not=\emptyset$$
When those two sets exist, the relation is written $\truthrelation_{\mathcal{D}_p, \mathcal{D}_c}$. 
\end{definition}

\begin{definition} A mixed truth-relation $\truthrelation_{\mathcal{D}_p, \mathcal{D}_c}$ is called \emph{$p$-mixed} iff $\mathcal{D}_{p} \subseteq \mathcal{D}_{c}$. 
\end{definition}

\begin{definition} A mixed truth-relation $\truthrelation_{\mathcal{D}_p, \mathcal{D}_c}$ is called \emph{$q$-mixed} iff $\mathcal{D}_{c} \subseteq \mathcal{D}_{p}$. 
\end{definition}

\begin{definition} A truth-relation $\truthrelation$ is called \emph{pure} if it is a mixed relation such that $\mathcal{D}_p=\mathcal{D}_c$, hence if it is both $p$-mixed and $q$-mixed. 
\end{definition}

\begin{definition}\label{def:sem2}
A semantics is called \emph{intersective mixed} (respectively \emph{intersective $p$-mixed} / \emph{intersective $q$-mixed} / \emph{intersective pure}) if it is truth-relational and its truth-relation is an intersection of mixed relations (respectively of $p$-mixed / $q$-mixed / pure relations).
\end{definition}

{A mixed semantics is trivially intersective mixed, but the converse is not true: an intersection of mixed semantics is not in general expressible as a mixed semantics, see \cite{CES2016} Theorem 2.15, with $ss\cap tt$ as a representative case in 3-valued semantics.} 

\subsubsection{A further semantic property: polarization}

Our definitions so far put no structure on the set of truth values. On an inferentialist perspective, however, it is often desirable to single out two special truth values, True and False, matching propositions with specific inferential roles. The value False is attached to the principle that from a contradictory proposition anything follows (\emph{ex falso quodlibet}), but also to the principle that if an argument is not valid, then the addition of a contradiction to the conclusions won't make it valid. The roles for the True are dual: a tautology should follow from any premise whatsoever, but if an argument is not valid, adding it among the premises will not make it valid either. 
Using $1$ and $0$ to represent the True and the False, those principles correspond to:

\begin{center}
\begin{tabular}{r@{~$\forall\gamma,\delta:$\hspace{3ex}}l@{\hspace{2em}}r@{~$\forall\gamma,\delta:$\hspace{3ex}}l}
(T1) & $\gamma\phantom{, 0} \truthrelation \delta, 1$ &
(T2) & $\gamma, 1 \truthrelation \delta\phantom{, 0}$ implies $\gamma \truthrelation \delta$ 
\\
(F1) & $\gamma, 0 \truthrelation \delta\phantom{, 1}$ &
(F2) & $\gamma\phantom{, 1} \truthrelation \delta, 0$ implies $\gamma \truthrelation \delta$ 
\end{tabular}
\end{center}

However, none of the semantic properties we have introduced so far for truth-relations guarantees that there will be truth values playing the roles of the True and the False. Consider for instance a relation of the form $\truthrelation_{\{\alpha\},\{\beta\}}$, requiring that when all premises take the value $\alpha$, some conclusion must take the distinct value $\beta$. The value $\beta$ is the only one that could follow from any set of values ($\truthrelation \beta$), but it is not the case that adding it as a premise has no effect on validity (e.g., $\alpha\not \truthrelation \emptyset$, but $\alpha,\beta \truthrelation \emptyset$). One can construct other odd truth-relations in which there is no room for special values such as {True and False}. The following property blocks some of these possibilities:

\begin{definition}[Polarization, True and False]
A mixed truth-relation $\truthrelation_{\mathcal{D}_p,\mathcal{D}_c}$ is \emph{polarized} iff it is \emph{T-polarized} and \emph{F-polarized}, in the following sense:
	\[
	\textrm{\textbf{T-polarization: }} \mathcal{D}_p\cap\mathcal{D}_c\not=\emptyset \quad\quad\quad
	\textrm{\textbf{F-polarization: }} \mathcal{V}\setminus(\mathcal{D}_p\cup\mathcal{D}_c)\not=\emptyset
	\]

The elements of $\mathcal{D}_p\cap\mathcal{D}_c$ are all candidates for {True}; the elements of $\mathcal{V}\setminus(\mathcal{D}_p\cup\mathcal{D}_c)$ are all candidates for {False}. An intersective mixed truth-relation is called \emph{polarized} if one of the mixed relations of which it is the intersection is T-polarized and one is F-polarized.%
\footnote{This is a weak requirement given the discussion above. It will be sufficient for all relevant purposes because most of the time we will be able to replace all intersective mixed consequence relations with mixed consequence relations.}
\end{definition}

The virtues of polarization become especially apparent when put in parallel with the structural property introduced in Definition~\ref{def:permeable}. There, we argued that permeable logics were odd, even if the standard Tarskian constraints on consequence relations do not rule them out. The following result shows that ruling out permeability at the structural level involves the admission of values playing the role of True and False at the semantic level.

\begin{theorem}\label{thm:permeablepolarized}
A sound and complete intersective mixed semantics for a non-permeable logic is polarized.
\end{theorem}

\begin{myproof}
Consider a sound and complete intersective mixed semantics for a non-permeable logic, {such that $\truthrelation=\bigcap_{\lambda \in \Lambda}\truthrelation_{\mathcal{D}_p^\lambda,\mathcal{D}_c^\lambda}$}. The failure of left-to-right permeability corresponds to F-polarization, and the failure of right-to-left permeability corresponds to T-polarization.
Choose $\Gamma, \Delta, \Sigma$ such that $\Gamma, \Sigma \vdash \Delta$ and $\Gamma \not\vdash \Sigma, \Delta$. Hence, there is some world $w$ such that:
	$\interp{\Gamma}(w)\not\truthrelation \interp{\Sigma, \Delta}(w)$
	and
	$\interp{\Gamma, \Sigma}(w)\truthrelation \interp{\Delta}(w)$. 
It follows that there is some member $\truthrelation_{\mathcal{D}_p^\lambda,\mathcal{D}_c^\lambda}$ of $\truthrelation$ such that:
	$\interp{\Sigma}(w)\cap\mathcal{D}_c^\lambda=\emptyset$ 
	and 
	$\interp{\Sigma}(w)\not\subseteq\mathcal{D}_p^\lambda$, 
i.e.~there is a truth value, in fact any of those in $\interp{\Sigma}(w)$, which is not in  $\mathcal{D}_p^\lambda$, and not in $\mathcal{D}_c^\lambda$ either. This proves F-polarization.
Similarly, choose $\Gamma, \Delta, \Sigma$ such that $\Gamma \vdash \Sigma, \Delta$ and $\Gamma, \Sigma \not\vdash \Delta$. Hence, there is some world $w$ such that:
	$\interp{\Gamma, \Sigma}(w)\not\truthrelation \interp{\Delta}(w)$ 
	and 
	$\interp{\Gamma}(w)\truthrelation \interp{\Sigma, \Delta}(w)$. 
It follows that there is some member $\truthrelation_{\mathcal{D}_p^\lambda,\mathcal{D}_c^\lambda}$ of $\truthrelation$ such that:
	$\interp{\Sigma}(w)\subseteq\mathcal{D}_p^\lambda$ 
	and 
	$\interp{\Sigma}(w)\cap\mathcal{D}_c^\lambda\not=\emptyset$, 
i.e.~$\mathcal{D}_p^\lambda\cap\mathcal{D}_c^\lambda\not=\emptyset$, which proves T-polarization.
\end{myproof}

This result is the first connection we establish between a structural property and a semantic property. In the next section, we proceed to flesh out more systematic parallelisms between semantic properties and structural properties, a summary of which is given in Table~\ref{tab:representation}.

\section{Correspondence between {structural} and semantic properties}\label{sec:correspondence}

In this section we state and prove representation theorems establishing a systematic correspondence between structural properties of a logic and fundamental properties of a semantics it may be associated with. 
We start with a semantic characterization of substitution-invariance. We find antecedent characterizations, in particular in the work of W\'ojcicki (see \citealp{bloom1970theorems, wojcicki1973matrix, wojcicki1988theory}), but in the context of Tarskian consequence relations (using logical matrices). Our result, on the other hand, offers a semantic counterpart of substitution-invariance, singled out from other Tarksian properties. The next result, Theorem \ref{thm:mon}, may be viewed as the cornerstone, for it shows that monotonicity, among all of Tarski's constraints, corresponds tightly to mixed consequence.
We then inspect further Tarskian properties, reflexivity, transitivity and their conjunction, which show up as particular types of mixed consequence. 
We give a summary of the main correspondences in Table \ref{tab:representation}. 

Closely related results were stated independently by \cite{french2017valuations} (Theorem~2), {and by \cite{blasio2017inferentially} ({Theorem 10}, Theorem 11)}, see also \cite{humberstone1988heterogeneous} for a precursor. 
Even though formulated and proved differently, {some of these results} are exactly equivalent in so far as substitution-invariance is dropped {(an exception is \citealt{blasio2017inferentially})} and, correspondingly, in so far as its semantic counterpart valuationality is dropped (replace strong semantics with truth-adequate semantics in our formulations; we thus give two versions of those results in what follows, versions marked with a prime drop the substitution-invariance assumption). Overall, they establish a similar correspondence between structural properties of consequence relations and {semantic properties}.

\begin{table}[!h]%
\[
\begin{tabular}{l|l}
\multicolumn{2}{c}{\textbf{Logic}}\\
\hline
\multicolumn{2}{c}{Substitution-invariant}\\
\multicolumn{2}{c}{Monotonic}\\
Reflexive & Transitive\\
\multicolumn{2}{c}{Tarskian}\\

\end{tabular}\qquad
\begin{tabular}{l|l}
\multicolumn{2}{c}{\textbf{Semantics}}\\
\hline
\multicolumn{2}{c}{Strong}\\
\multicolumn{2}{c}{Intersective Mixed}\\
Intersective $p$-Mixed & Intersective $q$-Mixed\\
\multicolumn{2}{c}{Intersective Pure}\\

\end{tabular}
\]

\caption{
Correspondence between structural properties of a logic and semantic properties of semantics for it. The properties should be read cumulatively from top to bottom, and existentially from left to right: substitution-invariant logic = existence of a strong semantics; (substitution-invariant + monotonic) logic = existence of a (strong + intersective mixed) semantics; etc. 
}
\label{tab:representation}
\end{table}

\subsection{Substitution-invariance = Strong semantics}

\begin{theorem}\label{th:structural=valuational}
A logic is substitution-invariant iff it has a strong semantics.
\end{theorem}

\begin{myproof}
\begin{description}
\item[($\Leftarrow$)]
	Assume that a logic $\langle\mathcal{L}, \mathcal{C}, \vdash\rangle$ has 
	a strong semantics $\langle\mathcal{V}, \interp{\_}\rangle$. 
	Suppose $\Gamma\vdash\Delta$, and $\sigma$ is a substitution. Then:

	$\forall v: v(\Gamma)\truthrelation v(\Delta)$ (by truth-relationality and adequacy), hence

	$\forall v: v\circ\sigma(\Gamma)\truthrelation v\circ\sigma(\Delta)$ (by valuationality: if $v$ is a valuation, then $v\circ\sigma$ is a valuation); moreover 
		
	$\forall v: v(\Gamma[\sigma])\truthrelation v(\Delta[\sigma])$ 
	($v\circ\sigma(A)=v(A[\sigma])$ for all formulae $A$, by valuationality for $A$ an atom, and by truth-functionality for the induction for more complex formulae), 
	
	so $\Gamma[\sigma] \vdash \Delta[\sigma]$ (by truth-relationality and adequacy).

\item[($\Rightarrow$)] Assume that a logic $\langle\mathcal{L}, \mathcal{C}, \vdash\rangle$ is subtitution-invariant. We define a strong semantics $\langle\mathcal{V}, \interp{\_}\rangle$ with 
	$\mathcal{V}=\mathcal{L}$, 
	$\sem{\vdash}={\vdash}$, 
	$\sem{c}=c$. 
	Valuationality, truth-interpretability and compositionality are built in (see Fact~\ref{fact:semanticsfromtruth}) by completing the semantics: we consider all valuations and assign an interpretation inductively for complex formulae. Let us prove soundness and completeness:
	$\Gamma\vdash\Delta$
	iff $\forall \sigma: \Gamma[\sigma]\vdash\Delta[\sigma]$ (by substitution-invariance),
	iff $\forall \sigma: \sigma(\Gamma)\truthrelation\sigma(\Delta)$ (by construction for all $F$: $\sigma(F)=F[\sigma]$ here), and substitutions {are} worlds/valuations (worlds here are all functions from $\mathcal{A}$ to the set of truth values, which is $\mathcal{L}$ here).
\qedhere
\end{description}
\end{myproof}

\subsection{Monotonic relations = intersection of mixed relations}

\begin{theorem}\label{thm:mon}
A logic is substitution-invariant and monotonic iff it has a strong intersective mixed semantics.
	\end{theorem}	

\begin{thmbis}{thm:mon}\label{thm:monNOSI}
A logic is monotonic iff it has a truth-adequate intersective mixed semantics.
\end{thmbis}	

\begin{myproof}
\begin{description}

\item [($\Leftarrow$)] Assume that we have a strong intersective mixed semantics.
Substitution-invariance follows from Theorem~\ref{th:structural=valuational}. The monotonicity of the truth-relation and therefore of the consequence relation is straightforward based on Definition \ref{def:mixed}.

\item [($\Rightarrow$)]

Consider a substitution-invariant and monotonic logic. We will construct a strong and intersective mixed semantics. For this, we use a method that closely follows methods attributed to Adolph Lindenbaum, used by, e.g., Tarski and W\'ojcicki and described in, e.g., \cite{bloom1970theorems, surma1982origin}.

Take $\mathcal{V}=\mathcal{L}$ and let the interpretation function $\interp{\_}$ on formulae and connectives be as in Theorem~\ref{th:structural=valuational} under ($\Rightarrow$).

Note $\lambda=(\Gamma, \Delta)$ such that $\Gamma\not\vdash\Delta$. Then define $\truthrelation_{\lambda}$ as the mixed consequence truth-relation on the sets: $\mathcal{D}_p^\lambda=\Gamma$ and $\mathcal{D}_c^\lambda=\mathcal{V}\setminus\Delta$. Define $\truthrelation$ as the intersection of these mixed consequence relations (to account for the universal consequence relation, i.e.~with no $\Gamma\not\vdash\Delta$: we will consider that it is indeed the intersection of an empty set of mixed consequence relations, alternatively any semantics with an empty set of worlds $\mathcal{W}$ would be appropriate). The resulting semantics is valuational, compositional, truth-interpretable and intersective mixed, by construction. We can show that it is sound and complete for the logic, in the sense that $\Gamma\vdash\Delta$ iff $\forall v, \forall \lambda: v(\Gamma)\truthrelation_\lambda v(\Delta)$.

\begin{itemize}
\item Assume $\Gamma\not\vdash\Delta$. Then take $\lambda=(\Gamma, \Delta)$, $v$ the identity valuation, it is clear that $v(\Gamma)\not\truthrelation_\lambda v(\Delta)$ and therefore that it's not the case that $\forall\lambda: \forall v:v(\Gamma)\truthrelation_\lambda v(\Delta)$.
\item Assume $\Gamma\vdash\Delta$. Then suppose it's not the case that $\forall\lambda', \forall v:v(\Gamma)\truthrelation_{\lambda'} v(\Delta)$. That is, there is $\lambda'=(\Gamma', \Delta')$ such that $\Gamma'\not\vdash\Delta'$ and $v$ such that $v(\Gamma)\subseteq \Gamma'$ and $v(\Delta)\subseteq \Delta'$. But $v(\Gamma)$ and $v(\Delta)$ are simply obtained by atomic substitutions in $\Gamma$ and $\Delta$, and therefore $v(\Gamma)\vdash v(\Delta)$ holds by substitution-invariance. By monotonicity we would then obtain that $\Gamma'\vdash\Delta'$, which was excluded.
\qedhere
\end{itemize}

\end{description}

\end{myproof}

\begin{myproofbis}is obtained similarly, by considering the same construction but reduced to a single valuation, the identity valuation, as in Fact~\ref{fact:t-adeqforall}.
\end{myproofbis}

	\subsection{Monotonic + reflexive relations = intersection of $p$-mixed relations}
		
	\begin{theorem}\label{thm:monref}
	
	A logic is substitution-invariant, monotonic and reflexive iff 
	it has a strong intersective $p$-mixed semantics.
	\end{theorem}
	
\begin{thmbis}{thm:monref}\label{thm:monrefNOSI}
A logic is monotonic and reflexive iff it has a truth-adequate intersective $p$-mixed semantics.
\end{thmbis}	

	\begin{myproof}
	($\Leftarrow$) Clearly every $p$-mixed consequence relation is reflexive, and so is an intersection of $p$-mixed relations.
($\Rightarrow$) That is because in the previous proof, we take 
$\truthrelation_{(\Gamma, \mathcal{V}\setminus\Delta)}$
for $\Gamma\not\vdash\Delta$, but if the relation is reflexive the latter requires that $\Gamma\cap\Delta=\emptyset$, therefore $\Gamma \subseteq \mathcal{V} \setminus \Delta$, hence $\mathcal{D}_p^\lambda\subseteq \mathcal{D}_c^\lambda$, making the relation $p$-mixed for every $\lambda$.
	\end{myproof}

\begin{myproofbis}is obtained similarly, by considering the same construction but reduced to a single valuation, the identity valuation, as in Fact~\ref{fact:t-adeqforall}.
\end{myproofbis}

	\subsection{Monotonic + transitive relations = intersection of $q$-mixed relations}

	\begin{theorem}\label{thm:montrans}
	
		A logic is substitution-invariant, monotonic and transitive iff 	
		it has a strong intersective $q$-mixed semantics.
	
	\end{theorem}	

\begin{thmbis}{thm:montrans}\label{thm:montransNOSI}
A logic is monotonic and transitive iff it has a truth-adequate intersective $q$-mixed semantics.
\end{thmbis}

	\begin{myproof}
	\begin{description}

	\item[($\Leftarrow$)] Consider a strong intersective $q$-mixed semantics, and $\Gamma\not\vdash\Delta$. Hence, there is $\lambda$ and $v$ such that $v(\Gamma)\subseteq\mathcal{D}^{\lambda}_p$ and $v(\Delta)\cap\mathcal{D}^{\lambda}_c=\emptyset$. Let $\Gamma'=\{P|v(P)\in\mathcal{D}^{\lambda}_p\}$. Then $\Gamma, \Gamma'\not\vdash \Delta, \mathcal{L}\setminus(\Gamma')$, which proves transitivity.
	The reason is that (i)~$v(\Gamma,\Gamma')\subseteq\mathcal{D}_p^{\lambda}$, by construction and (ii)~$v(\Delta,{\mathcal{L}\setminus(\Gamma')})\cap\mathcal{D}^{\lambda}_c=\emptyset$ because on the one hand, $v(\Delta)\cap\mathcal{D}^{\lambda}_c=\emptyset$ (see above), and on the other hand, $v(\mathcal{L}\setminus(\Gamma'))\cap\mathcal{D}^{\lambda}_c\subseteq v(\mathcal{L}\setminus(\Gamma'))\cap\mathcal{D}^{\lambda}_p=\emptyset$, {since $D_c^\lambda \subseteq D_p^\lambda$}. 

	\item[($\Rightarrow$)] Now consider a substitution-invariant, monotonic and transitive logic. We can construct a valuational and intersective $q$-mixed truth semantics. We will use the same type of `Lindenbaum' construction as used in the proof of Theorem~\ref{thm:mon}, but we will only consider a subset of the truth-relations for the intersection. 
	Consider the $\widetilde{\lambda}=(\widetilde{\Gamma},\widetilde{\Delta})$ such that $\widetilde{\Gamma}\not\vdash\widetilde{\Delta}$ \emph{but now also} $\widetilde{\Gamma}\cup\widetilde{\Delta}=\mathcal{V}$. We will consider $\truthrelation_{\widetilde{\Lambda}}$, the intersection of the $\truthrelation_{\widetilde{\lambda}}$s defined as the mixed relation with $\mathcal{D}_p^{\widetilde{\lambda}}=\widetilde{\Gamma}$ and $\mathcal{D}_c^{\widetilde{\lambda}}=\mathcal{V}\setminus\widetilde{\Delta}$ (as before).
	The rest of the proof amounts to showing soundness and completeness, it relies on three results:
	\begin{enumerate}
	\item The $\truthrelation_{\widetilde{\lambda}}$s are $q$-mixed: since $\widetilde{\Gamma}\cup\widetilde{\Delta}=\mathcal{V}$, $(\mathcal{V}\setminus\widetilde{\Delta})\subseteq\widetilde{\Gamma}$, i.e. $\mathcal{D}_c^{{\widetilde{\lambda}}}\subseteq\mathcal{D}_p^{{\widetilde{\lambda}}}$.
	\item If $\Gamma\vdash\Delta$ then $\forall v: v(\Gamma)\vdash v(\Delta)$: this holds for $\truthrelation_{\Lambda}$ (see above) and therefore also for $\truthrelation_{\widetilde{\Lambda}}$ which is weaker since it is an intersection of fewer truth-relations.
	\item If $\Gamma\not\vdash\Delta$, then 
	by transitivity we can find $\widetilde{\lambda}=(\widetilde{\Gamma},\widetilde{\Delta})$ such that 
$\widetilde{\Gamma}\not\vdash\widetilde{\Delta}$ and $\widetilde{\Gamma}\cup\widetilde{\Delta}=\mathcal{V}$. 
	$\Gamma\not\truthrelation_{(\widetilde{\Gamma},\mathcal{V}\setminus\widetilde{\Delta})}\Delta$: this is so because by construction $\widetilde{\Gamma}\not\truthrelation_{(\widetilde{\Gamma},\mathcal{V}\setminus\widetilde{\Delta})}\widetilde{\Delta}$ and $\Gamma\subseteq\widetilde{\Gamma}$ and $\Delta\subseteq\widetilde{\Delta}$.
	\qedhere
	\end{enumerate}
	\end{description}
	\end{myproof}

\begin{myproofbis}is obtained similarly, by considering the same construction but reduced to a single valuation, the identity valuation, as in Fact~\ref{fact:t-adeqforall}.
\end{myproofbis}

	\subsection{Monotonic + reflexive + transitive = intersection of pure relations} 	
	\begin{theorem}\label{thm:monrefltrans}
	A logic is substitution-invariant, monotonic, reflexive and transitive iff 
	it has a strong intersective pure semantics.
	
	\end{theorem}

\begin{thmbis}{thm:monrefltrans}\label{thm:monrefltransNOSI}
A logic is monotonic, reflexive and transitive iff it has a truth-adequate intersective pure semantics.
\end{thmbis}	
	
	\begin{myproof}
	($\Leftarrow$) follows from the previous two theorems. ($\Rightarrow$) Consider the construction from the proof of Theorem~\ref{thm:montrans}. The $\widetilde{\lambda}$ must be compatible with reflexivity, that is $\widetilde{\Gamma}\cap\widetilde{\Delta}=\emptyset$. Hence we obtain not only that $\widetilde{\Gamma}\subseteq(\mathcal{V}\setminus\widetilde{\Delta})$ as in Theorem~\ref{thm:montrans}, but also necessarily that $\widetilde{\Gamma}\supseteq(\mathcal{V}\setminus\widetilde{\Delta})$ as in Theorem~\ref{thm:monref}. 
	\end{myproof}

\begin{myproofbis}is obtained similarly, by considering the same construction but reduced to a single valuation, the identity valuation, as in Fact~\ref{fact:t-adeqforall}.
\end{myproofbis}

\section{Ranks and truth-functional reduction} 

What is the least number of values needed to semantically represent a logic? Suszko's response to the question was two. Malinowski answered three, but relaxing some of the structural constraints that Suszko took for granted for a consequence relation. On our view, the answer to this question depends not just on the structural desiderata that are put on a consequence relation, but also on the desiderata that are placed on the semantics. In this section, we start out by reviewing Suszko's proposal: we basically explain the sort of reduction that both Scott and Suszko used in order to reduce truth values to a minimum `logical' rank (\citealp{scott1974completeness, suszko1977fregean}). As observed by several authors and first by Suszko himself (\citealp{beziau2001sequents, suszko1977fregean, caleiro2003dyadic, marcos2009non, caleiro2012many, font2009taking, shramko2011truth}), this kind of reduction on the number of algebraic truth values meets weak demands, in particular regarding compositionality (see Appendix~\ref{app:SSredcompos} for a more nuanced view). Although the reduction does not respect truth-compositionality in general, we show that a modification of it does so for compact logics operating with what we call \emph{regular} connectives. 
In an Appendix, Appendix~\ref{app:groupingreductions}, we investigate another notion of reduction, namely \emph{grouping reduction}, which does preserve many of the relevant structural properties we are interested in in all cases, but also incurs limitations. {{Importantly, in what follows we focus mostly on notions of rank, namely the minimum number of truth values needed to represent a logic. We do not present actual Suszko reductions of specific logics for that matter, but mostly use the notion toward irreducibility arguments. In that regard, our perspective differs from the 
important work of \cite{caleiro2012many} and \cite{caleiro2015bivalent}, who propose a constructive method to compute Suszko reductions from arbitrary semantics. Their approach is in a sense more general, because it is not restricted to logics for which the reduction preserves truth-functionality, and it is output-oriented. Our goal is different, however, and input-oriented: it is to identity constraints on the connectives of the base logic that will deliver a truth-functional reduction. This is what leads us to isolate the notion of a regular logic.

}

\subsection{The Suszko rank and the Scott-Suszko reduction}

	\begin{definition}[Suszko rank]\label{def:Suszkorank}
	The \emph{Suszko rank} of a logic is the least cardinal of truth values needed to get a sound and complete intersective mixed semantics for that logic.%
	\end{definition}

The restriction to intersective mixed semantics was not explicitly part of the issue in some past descriptions. But it is unproblematic given the results above for monotonic logics. Furthermore, previous results were always stated through \emph{mixed} semantics, and we therefore propose immediately a stronger notion of rank, which will prove not to be very different:

	\begin{definition}[mixed Suszko rank]\label{def:mSuszkorank}
	The \emph{mixed Suszko rank} of a logic is the least cardinal of truth values needed to get a sound and complete mixed semantics for that logic.
	\end{definition}

	Scott and Suszko came up with a method which outputs ranks 2, 3 or 4, depending on the structural constraints assumed to begin with. In their original approach, restricted to rank 2, this method consists in defining what Suszko calls \emph{logical valuations} (and Scott \emph{truth valuations}). Logical valuations are functions from formulae to the values $\{1,0\}$, defined on the basis of standard \emph{algebraic valuations}, namely functions from formulae to algebraic truth values. A logical valuation is derived from an algebraic valuation $h$ and a pure truth-relation, that is a set of designated algebraic values: it assigns a formula $F$ the value $1$ if $h$ assigns $F$ a designated value, and $0$ otherwise. In our framework, the reduction can be described as follows:
		
	\begin{definition}[Scott-Suszko reduction]\label{def:Scott-Suszko-reduction}
	Consider some intersective mixed semantics $\langle\mathcal{V}, \mathcal{W}, \interp{\_}\rangle$ for a logic, assuming that the truth-relation is given by the intersection of the $\truthrelation_{\mathcal{D}_p^{\lambda}, \mathcal{D}_c^{\lambda}}$ for $\lambda=(\mathcal{D}_p^{\lambda}, \mathcal{D}_c^{\lambda})$ in $\Lambda$. 
	Its associated \emph{Scott-Suszko reduction} is the semantics defined as follows, where we leave aside the interpretation of the connectives (which may be interpreted arbitrarily for now).

	$\mathcal{V}^*:=\{1,\#_p,\#_c,0\}$, 
	$\quad$
	$\mathcal{W}^*:=\Lambda \times \mathcal{W}$,
	$\quad$	
	$\truthrelation^*:=\truthrelation_{\{1,\#_p\},\{1,\#_c\}}$,
	
	$\interp{F}^*(\lambda, v):=t_\lambda(\interp{F}(v))$, with 
	\begin{tabular}[t]{llll}
	$t_\lambda:$ & $\mathcal{V}$ & $\longrightarrow$ & $\mathcal{V}^*$\\
	& $\alpha$ & $\longmapsto$ & 
	{\scriptsize$\left\{\begin{array}{ll}
	1 & \textrm{if } \alpha\in\mathcal{D}_p^{\lambda}\cap\mathcal{D}_c^{\lambda} \\
	0 & \textrm{if }  \alpha\not\in\mathcal{D}_p^{\lambda}\cup\mathcal{D}_c^{\lambda} \\
	\#_p & \textrm{if } \alpha\in\mathcal{D}_p^{\lambda}\setminus\mathcal{D}_c^{\lambda} \\
	\#_c & \textrm{if } \alpha\in\mathcal{D}_c^{\lambda}\setminus\mathcal{D}_p^{\lambda} 
	\end{array}\right.$}
	\end{tabular}
	
	\end{definition}

	\begin{theorem}\label{thm:suszkoscott}
	The Scott-Suszko reduction of a sound and complete, intersective mixed semantics is a sound and complete, mixed semantics.
	\end{theorem}
	
	\begin{myproof} 
	This follows from the examination of the two ends in the following chain of equivalence:
	
	$\phantom{\Leftrightarrow:}$ 
	$\forall (\lambda, v)\in \mathcal{W}^{*}: \interp{\Gamma}^{*}(\lambda, v)\truthrelation_{\{1,\#_p\},\{1,\#_c\}} \interp{\Delta}^{*}(\lambda, v)$
	
	$\Leftrightarrow:$
	$\forall (\lambda, v)\in \mathcal{W}^{*}: t_\lambda(\interp{\Gamma}(v))\truthrelation_{\{1,\#_p\},\{1,\#_c\}}t_\lambda(\interp{\Delta}(v))$
	
	$\Leftrightarrow: \forall (\lambda, v)\in \mathcal{W}^{*}:$ 
	if $t_\lambda(\interp{\Gamma}(v))\subseteq \{1,\#_p\}$, 
	then $t_\lambda(\interp{\Delta}(v))\cap \{1,\#_c\}\neq \emptyset$
	
	$\Leftrightarrow: \forall (\lambda, v)\in \mathcal{W}^{*}:$ if $\interp{\Gamma}(v)\subseteq \mathcal{D}_p^{\lambda}$ then $\interp{\Delta}({v})\cap \mathcal{D}_c^{\lambda}\neq \emptyset$
	
	$\Leftrightarrow: \forall v: \interp{\Gamma}(v) (\bigcap_{\Lambda}\truthrelation_{\mathcal{D}_p^{\lambda}, \mathcal{D}_c^{\lambda}})\interp{\Delta}(v)$
	
	$\Leftrightarrow: \Gamma \vdash \Delta$ (in virtue of the truth-relationality and soundness and completeness of the input semantics).
	\end{myproof}

	The following corollary puts upper bounds on Suszko rank for monotonic consequence relations with various syntactic properties. This is the exact analog of Corollary 1 in \cite{french2017valuations}, {and of statements (M1) and (M2) in \cite{blasio2017inferentially}}:
	
	\begin{corollary}[Maximum Suszko rank theorems]\label{corol:maxSrank}
	A monotonic consequence relation is of Suszko rank:
	\begin{itemize}
	\item at most 2 if it is reflexive and transitive,
	\item at most 3 if it is transitive or reflexive,
	\item at most 4 in general.
	\end{itemize}
	\end{corollary}
	
	\begin{myproof}
	
	The $t_\lambda$s from the Scott-Suszko reduction may yield at most 4 values. In fact, 3 values are sufficient if the relation is reflexive (one can choose an intersective semantics with, for all $\lambda$s, $\mathcal{D}_p^{\lambda}\setminus\mathcal{D}_c^{\lambda}=\emptyset$, by Theorem~\ref{thm:monrefNOSI}) or transitive (one can choose $\mathcal{D}_c^{\lambda}\setminus\mathcal{D}_p^{\lambda}=\emptyset$, by Theorem~\ref{thm:montransNOSI}). Finally, 2 values are sufficient if the relation is transitive and reflexive because we can then choose a representation such that $\mathcal{D}_p^{\lambda}=\mathcal{D}_c^{\lambda}$ (Theorem~\ref{thm:monrefltransNOSI}).
	\end{myproof}

	We can further show that these upper bounds are reached. Malinowski famously did it for rank 3 and transitive relations. {The following gives an example of a 4-valued logic which is not further reducible.}

	\begin{fact}\label{fact:reaching4}
	There exist substitution-invariant monotonic logics of Suszko rank exactly 4.
	\end{fact}
	
	\begin{myproof}
	Consider the language of propositional logic, the valuational semantics $\langle\{0,1,\#_p,\#_c\},\interp{\_}\rangle$ such that $\truthrelation=\truthrelation_{\{1,\#_p\},\{1,\#_c\}}$. We may give an interpretation to connectives (e.g., conjunction and disjunction interpreted as $\min$ and $\max$ with $0$ being lower than any other value, $1$ being higher, and by definition $\min(\#_p,\#_c)=0$, $\max(\#_p,\#_c)=1$, negation could be defined as, e.g., $\neg \#_p=\neg \#_c=0$). But most importantly, we add to the language four $0$-ary connectives, which will receive constant truth-functions in this semantics. {This semantics can be used to complement the language with a consequence relation} 
	The resulting logic cannot be given a truth-relational semantics with less than 4 semantic values.
	
	Pick an intersective semantics. Given the syntactic properties of the formulae $\top$, $\bot$, $\#_p$ and $\#_c$, there are constraints on the truth values they may receive. For instance, $\vdash\top$ requires that all the truth values received by $\top$ be in all of the $\mathcal{D}_c^{\lambda}$s, and the fact that $\top\not\vdash$ requires that at least one truth value received by $\top$ be outside of at least one $\mathcal{D}_c^{\lambda_0}$. Further similar constraints are given in the table below; each line corresponds to the constraints imposed on the truth values of the formula given in the first column, as derived from each of the statements between parentheses.
\[
\framebox{\begin{tabular}{l|l@{: }ll@{: }l}
$\top$ 
	&  ($\emptyset\vdash\top$) & all values in all $\mathcal{D}_c$
	& ($\top\not\vdash\emptyset$) & some value outside some $\mathcal{D}_p$ \\

$\bot$ 
	& ($\bot\vdash$) & all values outside all $\mathcal{D}_p$ 
	& ($\not\vdash\bot$) & some value in some $\mathcal{D}_c$  \\

$\#_c$ 
	& ($\#_c\vdash$) & all values outside all $\mathcal{D}_p$
	& ($\vdash\#_c$) & all values in all $\mathcal{D}_c$ \\

$\#_p$ 
	& ($\#_p\not\vdash\#_p$) & \multicolumn{3}{l}{some value in a given $\mathcal{D}_p$ but not in its $\mathcal{D}_c$ counterpart} \\
	
\end{tabular}}\]

It follows that the semantics has at least four values. First, it has a value $x_1$ in all $\mathcal{D}_c$s but outside some $\mathcal{D}_p$ (first line). Second, it has a value $x_2$ outside all $\mathcal{D}_p$s but in some $\mathcal{D}_c$, this cannot be $x_1$ which is in some $\mathcal{D}_p$. Third, it has a value $x_3$ outside all $\mathcal{D}_p$s (hence it is not $x_1$) and outside all $\mathcal{D}_c$s (hence it is not $x_2$). Finally, it has a value $x_4$ in a given $\mathcal{D}_p$ but not in its $\mathcal{D}_c$ counterpart. This cannot be $x_1$ because $x_1$ is in all $\mathcal{D}_c$s, it cannot be $x_2$ or $x_3$  because they are both outside all $\mathcal{D}_p$s.
\qedhere
\end{myproof}

More {generally}, we can obtain an exact correspondence between mixed Suszko rank and structural properties for non-permeable logics:

	\begin{corollary}[Exact Suszko rank theorems]\label{corol:exactSrank}
	A monotonic, non-permeable consequence relation is of mixed Suszko rank:
	\begin{itemize}
	\item exactly 2 if and only if it is reflexive and transitive,
	\item exactly 3 if and only if it is transitive or reflexive but not both,
	\item exactly 4 if and only if it is neither transitive nor reflexive.
	\end{itemize}
	\end{corollary}

\begin{myproof}
Consider a truth-adequate mixed semantics of minimal rank for such a logic, we note its truth-relation $\truthrelation_{\mathcal{D}_p, \mathcal{D}_c}$. In virtue of Theorem~\ref{thm:permeablepolarized}, some truth value is in both $\mathcal{D}_p$ and $\mathcal{D}_c$ (call one of them $1$), and some value is in neither (call one of them $0$). 

\begin{itemize} 

\item 
The rank for all logics we are concerned with here is at least 2 (the set of truth values contains $0$ and $1$), and hence it is exactly 2 for reflexive and transitive logics (for which the rank is at most 2, by Corollary~\ref{corol:maxSrank}). 

\item If the semantics is two-valued, then the truth-relation is essentially of the form $\truthrelation_{\{1\},\{1\}}$. Hence, the truth-relation is pure and by Theorem~\ref{thm:monrefltransNOSI}, the logic is reflexive and transitive.

\item If the semantics is 3-valued, call $\#$ the third truth value. In virtue of Theorem~\ref{thm:permeablepolarized}, the truth-relation will be either: $\truthrelation_{\{1\},\{1\}}$, $\truthrelation_{\{1,\#\},\{1\}}$, $\truthrelation_{\{1\},\{1,\#\}}$ or $\truthrelation_{\{1,\#\},\{1,\#\}}$. It follows (Theorems~\ref{thm:monrefNOSI} and~\ref{thm:montransNOSI}) that the logic is reflexive or transitive (or both). Hence, logics which are neither transitive nor reflexive will be of rank at least 4, and given Corollary~\ref{corol:maxSrank}, of rank exactly~4.
\qedhere\end{itemize}
\end{myproof}

\subsection{Stronger notions of rank: the problem of {truth-functionality}}\label{sec:compositionalityproblem}

Suszko's focus was on the notion of a logical truth value, and only concerns the existence of a sound and complete semantics for a logic.
However, one may demand more from a semantics, such as valuationality, compositionality and truth-functionality {(see \citealt{font2009taking}, \citealt{marcos2009non})}. Such properties may be seen as irrelevant in relation to logical truth values, insofar as they concern connectives more than consequence relations, but they surely are relevant in relation to the notion of an algebraic truth value. The Suzsko problem is about finding the \emph{rank} of a logic, for more or less stringent notions of rank, depending on the demands of the semantics.

\begin{definition}[Notions of rank] We introduce a sample of possible definitions of logical/algebraic rank, i.e.~possible alternatives to Suszko rank from Definition~\ref{def:Suszkorank}:
\begin{itemize}
\item The \emph{compositional rank} of a logic is the least cardinal of truth values needed to construct an adequate intersective mixed semantics for it. 
\item The \emph{truth-{functional} rank} of a logic is the least cardinal of truth values needed to construct a truth-adequate intersective mixed semantics for it.
\item The \emph{strong rank} of a logic is the least cardinal of truth values needed to construct a strong intersective mixed semantics for it.
\end{itemize}
As we did for the Suszko rank and the mixed Suszko rank, we can define parallel and in fact very similar notions of \emph{mixed compositional rank}, \emph{mixed truth-
{functional} rank} and \emph{mixed strong rank}, by requiring in each of the above definitions that the semantics be mixed (and not only intersective mixed).
\end{definition}

Compositionality in our general sense is a very weak demand and a very minor modification of Scott and Suszko's reduction guarantees it, in a trivial sense. We make this explicit in Appendix~\ref{app:SSredcompos}, where the reader can also find a more direct construction to obtain a Scott-Suszko like semantics for a logic, without going through the trouble of the reduction of some pre-existing semantics (this is to be compared to \citeauthor{french2017valuations}'s approach).

Hence, in the remainder of this paper, we will be interested in truth-functionality, and therefore in the truth-{functional} rank, not in the compositional rank. In general, the truth-{functional} rank is strictly more demanding than the Suszko rank for monotonic logics: even though the Suszko rank is at most 2 for monotonic, reflexive and transitive logics (Corollary~\ref{corol:maxSrank}), for a variety of such logics the truth-{functional} rank is higher than 2, {as \cite{suszko1977fregean} emphasized, a point we illustrate with the following Fact}:

	\begin{fact}\label{thm:notrured}
	Some monotonic, reflexive and transitive logic has no bivalent, truth-adequate intersective mixed semantics. 
		\end{fact}

	\begin{myproof}
	Consider the language of propositional logic, with a trivalent semantics for the connectives and formulae (over the set $\{1, \frac{1}{2}, 0\}$), and the order-theoretic consequence relation, which here is simply based on the truth-relation $\truthrelation_{\{1\},\{1\}}\cap\truthrelation_{\{1, \frac{1}{2}\},\{1, \frac{1}{2}\}}$. 
	Also assume that it has a formula $\top$ for the tautology (constant proposition of value $1$) and one formula $\bot$ for the contradiction (constant proposition of value $0$). It is substitution-invariant,  monotonic, reflexive and transitive, since it has a valuational, intersective pure semantics. Consider, for reductio, a bivalent truth-adequate semantics. 
	
We first need to show that the truth-relation necessarily is the classical bivalent consequence truth-relation. 
The truth-relation $\truthrelation$ is intersective: it is the intersection of all $\truthrelation_{\mathcal{D}_p^{\lambda},\mathcal{D}_c^{\lambda}}$ such that $\mathcal{D}_p^{\lambda}\not\truthrelation(\mathcal{V}\setminus\mathcal{D}_c^{\lambda})$ (the reasoning is similar to the one in the proof of Theorem~\ref{thm:mon}).
	One of the two values in $\mathcal{V}$, call it $1$, has to be in all $\mathcal{D}_c^{\lambda}$s, for otherwise there could be no validity ($\vdash\top$) and in some $\mathcal{D}_p^{\lambda}$ (otherwise it cannot be that $\top\not\vdash$). One value, call it $0$, has to be in no $\mathcal{D}_p^{\lambda}$, for otherwise there could be no proposition that validates the empty set of conclusions ($\bot\vdash$) and outside of some $\mathcal{D}_c^{\lambda}$ ($\bot\vdash$). The truth-relation is therefore the intersection of some of the following:
	$\truthrelation_{\emptyset,\{1,0\}}$, $\truthrelation_{\emptyset,\{1\}}$, $\truthrelation_{\{1\},\{1,0\}}$, $\truthrelation_{\{1\},\{1\}}$. The last one is more stringent than any of the others, and it is the only one that allows for $\top\not\vdash\bot$ (the other relations hold as soon as the premises and the conclusions are not empty). Hence, $\truthrelation$ is $\truthrelation_{\{1\},\{1\}}$.
	
Then let us consider what the truth-function for negation would be. For an atom $p$: $p\not\vdash \neg p$. Hence, there is a valuation in which $p$ takes the value $1$ and $\neg p$ the value $0$, i.e., $\sem{\neg}(1)=0$. Similarly, $\neg p\not\vdash p$, from which it follows that $\sem{\neg}(0)=1$. But under those assumptions, we should have $\vdash p, \neg p$, which does not hold. Hence, no such semantics can be truth-functional.

	\end{myproof}

In fact, the discrepancy is more extreme, since even though the Suszko rank is at most 4 for any monotonic logics (Corollary~\ref{corol:maxSrank}), the strong rank is unbounded:

\begin{theorem}\label{thm:notrured-n}
For any $n$, there is some monotonic 
logic of truth-{functional} rank $n$, i.e.~with no truth-adequate semantics with fewer than $n$ truth values.
\end{theorem}

\begin{myproof}
Consider the language of propositional logic over a fixed set of connectives, with a semantics over $n$ totally ordered truth values $\{\alpha_1, ..., \alpha_n\}$, some truth-functions for the connectives, and the associated order-theoretic consequence relation, defined as $\gamma\truthrelation\delta$ iff $\min(\gamma)\leq\max(\delta)$. This truth-relation is also the intersection of all pure $\truthrelation_{\mathcal{D}_{i}, \mathcal{D}_{i}}$, with $\mathcal{D}_i=\{\alpha_k: k\geq i\}$, for $i\in\{2, ..., n\}$.
Also extend the language (and the semantics) with $n$ $0$-ary truth-functional connectives $C_1, ..., C_n$, with values corresponding to each of the $n$ truth values. This yields a logic which is substitution-invariant, monotonic, reflexive and transitive, since it has a valuational, intersective pure semantics.
In any truth-adequate semantics, the $0$-ary connectives correspond to a truth value (a constant truth-function), and they all have to be different truth values because the formulae $C_i$s have pairwise different behaviors ($C_i\vdash C_j$ iff $i\leq j$). Hence, any truth-adequate semantics has at least $n$ truth values.
\end{myproof}

Facts \ref{thm:notrured} and Theorem \ref{thm:notrured-n}  may be related to \cite{godel1932intuitionistic}'s result showing that there is no finite-valued characterization of intuitionist logic. Intuitionist consequence is monotonic, reflexive, transitive, and nonpermeable.\footnote{{For the intuitionistic calculus, \cite{gentzen1935investigations} allows only one formula to appear as conclusion in a sequent. This may appear to restrict monotonicity to premises only. However, Gentzen's restriction is primarily concerned with the application of operational rules, namely the rules for the connectives. In other words, an intuitionistically valid sequent can be augmented with more conclusions and remain valid under the intuitionistic reading of the comma as a disjunction on the conclusion side.
}}} From Corollary \ref{corol:exactSrank}, it follows that it must have a two-valued semantic characterization. There is no contradiction with G\"odel's result, however, for what G\"odel's result really establishes is that there is no \emph{truth-functional} characterization of intuitionistic logic with finitely many values.

\subsection{A solution for `regular' connectives and compact logics: the {truth-functional} Scott-Suszko reduction}

The counterexample given in the previous section to the truth-{functionality} of Scott-Suszko reductions involved connectives with particular properties: 
a negation for which it could be that $\Gamma, P \vdash \Delta$ without $\Gamma \vdash \neg P, \Delta$  (Theorem~\ref{thm:notrured})
or $0$-ary connectives specifically targeting algebraic truth-values, so to speak (Theorem~\ref{thm:notrured-n}).
Here we show that under certain conditions, the truth-{functional} rank is not higher than the Suszko rank. We exhibit a modification of the Scott-Suszko reduction that preserves truth-relationality and soundness and completeness, but also truth-{functionality}. But this happens only for compact logics and certain types of connectives that we will call regular. 

\subsubsection{Compact and regular logics}
Compactness is a standard structural notion from the Tarskian tradition:

\begin{definition}[Compactness]
A logic is \emph{compact} if $\Gamma\vdash\Delta$ entails that there are $\Gamma'$ and $\Delta'$ finite subsets of $\Gamma$ and $\Delta$, such that $\Gamma'\vdash\Delta'$.
\end{definition}

As for regular connectives, a formal definition is given in~\ref{def:regconn}. In essence, they are connectives such that their structural behavior can be read off from the consequence relation. The formulation of the definition is a bit complicated, but the spirit is simple and a couple of examples could help clarify. Whether a relation holds when the premise contains a formula headed by negation can, in well-behaved cases, be deduced from whether related relations without that negation hold. And similarly when the formula headed by negation is in conclusion. More generally, whether $\Gamma, C(F_1,..., F_n)\vdash \Delta$ is often dictated by a conjunction of statements of the form ``$\Gamma, F_1 \vdash F_2, \Delta$'', ``$\Gamma, F_2, F_3 \vdash F_1, \Delta$'', etc.
	
	\begin{example}\label{def:structconn}
	Examples of typical relations that make connectives \emph{regular}:\medskip
	
	\noindent \framebox{\begin{tabular}{l|@{ $\Gamma, $ }c@{ $\vdash \Delta\,$ iff $\,$}l|@{ $\Gamma \vdash$ }c@{ $, \Delta \,$ iff $\,$}l}
	Negation
		& $\neg P$ & $\Gamma \vdash P, \Delta$
		& $\neg P$ & $\Gamma, P \vdash \Delta$\\
	Conjunction
		& $P\wedge Q$ & $\Gamma, P, Q \vdash \Delta$
		& $P\wedge Q$ & $\Gamma \vdash P, \Delta$ and $\Gamma \vdash Q, \Delta$\\
	Disjunction
		& $P\vee Q$ & $\Gamma, P \vdash \Delta$ and $\Gamma, Q \vdash \Delta$
		& $P\vee Q$ & $\Gamma \vdash P, Q, \Delta$\\
	Conditional
		& $P{\to}{Q}$ &  $\Gamma, Q \vdash \Delta$ and $\Gamma \vdash P, \Delta$ 
		& $P{\to}{Q}$ & $\Gamma, P \vdash Q, \Delta$\\
	\end{tabular}}
	\end{example}

\noindent
All of those correspond to invertible sequent calculus rules for classical propositional logic (viz. \citealp{scott1974completeness}). The general definition of regular connectives is as follows.

	\begin{definition}[Regular connectives, regular logics]\label{def:regconn}

	A connective $C$ is \emph{regular} if there exist
	$\mathcal{B}^p\subseteq\mathcal{P}(\{1,..., n\})\times\mathcal{P}(\{1,..., n\})$ and
	$\mathcal{B}^c\subseteq\mathcal{P}(\{1,..., n\})\times\mathcal{P}(\{1,..., n\})$
such that
	$\forall\Gamma, \Delta:$
	\[\begin{array}{c@{\textrm{ iff }}c}
	\Gamma, C(F_1, ..., F_n) \vdash \Delta
		& \bigwedge\limits_{(B_p,B_c)\in \mathcal{B}^p} 
			{\Gamma, \{F_i: i\in B_p\}\vdash \{F_i: i\in B_c\}, \Delta}\\

	\Gamma \vdash C(F_1, ..., F_n) , \Delta
		& \bigwedge\limits_{(B_p,B_c)\in \mathcal{B}^c} 
			{\Gamma, \{F_i: i\in B_p\}\vdash \{F_i: i\in B_c\}, \Delta}\\
	\end{array}\]

	A logic is called \emph{regular} if all its connectives are regular. For a regular connective $C$, we may refer to rules as above as its \emph{regularity rules}, and note the corresponding sets $\mathcal{B}^p$ and $\mathcal{B}^c$, leaving $C$ implicit. We do so without implying that these are uniquely determined.
	\end{definition}

{Regular logics include not only classical logic, but also some substructural logics, such as the  logics ST and TS of \cite{tcs}.\footnote{{We refer to \cite{chemla2019mv} for a detailed investigation of regular connectives in such logics.}} Non-regular logics include for example intuitionistic logic or \L ukasiewicz's three-valued logic $\L 3$: in particular, in both systems the negation is not regular.\footnote{{The classic rule of Example \ref{def:structconn} fails in either logic. Moreover, none of the other options to make negation regular is valid: neither $\Gamma, \neg A\vdash \Delta$ iff $\Gamma\vdash A, \Delta$; nor $\Gamma, \neg A\vdash \Delta$ iff $\Gamma, A\vdash \Delta$; nor $\Gamma, \neg A\vdash \Delta$ iff $\Gamma\vdash \Delta$.}} Intuitively, the notion of regularity is tied to the proof-theoretical notion of analyticity, namely for a logic to enjoy the subformula property (all formulas occurring in a derivation of a provable sequent are subformulae of formulae in the end-sequent). We may conjecture that regularity is a sufficient condition for a logic to enjoy the subformula property. We leave a further exploration of that issue for another occasion (see \citealt{lahav2018subformula} for a recent study).}

\subsubsection{The {truth-functional} Scott-Suszko reduction and the {truth-functional} rank}

We can exploit the regularity of a logic to propose a version of the Scott-Suszko reduction designed to preserve {truth-functionality}, modifying the construction from Definition~\ref{def:Scott-Suszko-reduction}.
{In this reduction, we now say explicitly how to treat (regular) connectives.}

\begin{definition}[{Truth-functional} Scott-Suszko reduction]\label{def:strongssreduction}
	Consider some intersective mixed semantics $\langle\mathcal{V}, \mathcal{W}, \interp{\_}\rangle$ for a logic, assuming that the truth-relation is given by the intersection of the $\truthrelation_{\mathcal{D}_p^{\lambda}, \mathcal{D}_c^{\lambda}}$ for $\lambda=(\mathcal{D}_p^{\lambda}, \mathcal{D}_c^{\lambda})$ in $\Lambda$. Its \emph{{truth-functional} Scott-Suszko reduction} is the semantics defined as follows:

	$\mathcal{V}^*:=\{1,\#_p,\#_c,0\}$, 
	
	$\mathcal{W}^*:=\Lambda \times \mathcal{W}$,
	
	$\truthrelation^*\ :=\ \truthrelation_{\{1,\#_p\},\{1,\#_c\}}$,
	
	$\interp{p}^*(\lambda, v):=t_\lambda(\interp{p}(v))$ for all atoms $p$, and we extend the interpretation function inductively through the definition of a truth-function for all (regular) connectives:
	\[
	\begin{array}{c@{\textrm{ iff }}c}
	\sem{C}(x_1, ..., x_n)\not\in\{1,\#_p\} 
		& \bigwedge\limits_{(B_p,B_c)\in \mathcal{B}^p} 
			{(\{x_i: i\in B_p\} \subseteq\{1,\#_p\} \Rightarrow \{x_i: i\in B_c\} \cap\{1,\#_c\}\not=\emptyset)}\\

	\sem{C}(x_1, ..., x_n)\in\{1,\#_c\} 
		& \bigwedge\limits_{(B_p,B_c)\in \mathcal{B}^c} 
			{(\{x_i: i\in B_p\} \subseteq\{1,\#_p\} \Rightarrow \{x_i: i\in B_c\} \cap\{1,\#_c\}\not=\emptyset)}\\
	\end{array}\]

\end{definition}

This reduction relies on a canonical way of associating a truth-function to a regular connective. In Appendix~\ref{app:regconnectives}, we  provide reasons to compare the resulting truth-functions with Strong Kleene truth-tables. For current purposes, we want to ensure that these truth-functions deliver {truth-functionality}:

\begin{theorem}\label{thm:semredstructconn}
The {truth-functional} Scott-Suszko reduction of a sound and complete, intersective semantics for a regular and compact logic is a truth-adequate, mixed semantics. 
\end{theorem}
	
\begin{myproof}
By construction, the {truth-functional} Scott-Suszko reduction is truth-compositional and mixed, see Fact~\ref{fact:semanticsfromtruth}. We will now prove that it is sound and complete in three steps: the induced syntactic consequence relation 
(i)~is sound and complete on sets of premises and conclusions only made of atoms, 
(ii)~respects the regularity of the connectives, 
(iii)~is inductively sound and complete everywhere since by compactness and regularity all relations can be inferred from relations between atomic formulae only.
\begin{itemize}

\item[(i)] On atoms, the semantics is clearly sound and complete: it coincides with the traditional Scott-Suszko reduction.

\item[(ii)] This semantics respects the regularity of the connectives. We first show this for cases where the connective is in the premise set: the first member of the chain of equivalence below corresponds to 
	$\Gamma, C(F_1, ..., F_n) \vdash \Delta$ and the last member corresponds to its regular representation 
	$\bigwedge_{(B_p,B_c)\in \mathcal{B}^p} 
			{\Gamma, \{F_i: i\in B_p\}\vdash \{F_i: i\in B_c\}, \Delta}$. 
	$\forall v^*=(\lambda,v):$

{\scriptsize 	 
		$[ v^*(\Gamma)\subseteq\{1,\#_p\} 
		\wedge
		v^*(C(F_1, ..., F_n))\in\{1,\#_p\}]
		\Rightarrow
		v^*(\Delta)\cap\{1,\#_c\}\not=\emptyset$

	iff 
		$v^*(\Gamma)\not\subseteq\{1,\#_p\} 
		\vee
		v^*(C(F_1, ..., F_n))\not\in\{1,\#_p\}
		\vee
		v^*(\Delta)\cap\{1,\#_c\}\not=\emptyset$
	
	iff  
		$v^*(\Gamma)\not\subseteq\{1,\#_p\} 
		\vee
		v^*(\Delta)\cap\{1,\#_c\}\not=\emptyset
		\vee
		\bigwedge_{(B_p,B_c)\in \mathcal{B}^p}  
			{\{v^*(F_i): i\in B_p\} \subseteq\{1,\#_p\} 
			\Rightarrow 
			\{v^*(F_i): i\in B_c\} \cap\{1,\#_c\}\not=\emptyset}$

	iff  
		$\bigwedge_{(B_p,B_c)}
		[
		{\{v^*(F_i): i\in B_p\} \not\subseteq\{1,\#_p\} 
				\vee
		\{v^*(F_i): i\in B_c\} \cap\{1,\#_c\}\not=\emptyset}
		\vee
		v^*(\Gamma)\not\subseteq\{1,\#_p\} 
		\vee
		v^*(\Delta)\cap\{1,\#_c\}\not=\emptyset
		]
		$

	iff 
		$\bigwedge_{(B_p,B_c)}  
		[
		{\{v^*(F_i): i\in B_p\} \subseteq\{1,\#_p\} 
		\wedge
		v^*(\Gamma)\subseteq\{1,\#_p\} 
		]
		\Rightarrow
		[
		\{v^*(F_i): i\in B_c\} \cap\{1,\#_c\}\not=\emptyset}
		\textrm{ or }
		v^*(\Delta)\cap\{1,\#_c\}\not=\emptyset
		]
		$

}

Similarly, for cases where the connective is in the conclusion set: the first member of the chain of equivalence below corresponds to 
$\Gamma \vdash C(F_1, ..., F_n), \Delta$ and the last one to its regular representation 
	$\bigwedge_{(B_p,B_c)\in \mathcal{B}^c} 
		{\Gamma, \{F_i: i\in B_p\}\vdash \{F_i: i\in B_c\}, \Delta}$.
	$\forall v^*=(\lambda,v):$

{\scriptsize 	 
		$v^*(\Gamma)\subseteq\{1,\#_p\} 
		\Rightarrow
		[v^*(\Delta)\cap\{1,\#_c\}\not=\emptyset
		\vee
		v^*(C(F_1, ..., F_n))\in\{1,\#_c\}]$

	iff  
		$v^*(\Gamma)\not\subseteq\{1,\#_p\} 
		\vee
		v^*(\Delta)\cap\{1,\#_c\}\not=\emptyset
		\vee
		v^*(C(F_1, ..., F_n))\in\{1,\#_c\}$
	
	iff 
		$v^*(\Gamma)\not\subseteq\{1,\#_p\} 
		\vee
		v^*(\Delta)\cap\{1,\#_c\}\not=\emptyset
		\vee
		\bigwedge_{(B_p,B_c)\in \mathcal{B}^c}  
			{\{v^*(F_i): i\in B_p\} \subseteq\{1,\#_p\} \Rightarrow \{v^*(F_i): i\in B_c\} \cap\{1,\#_c\}\not=\emptyset}$

	iff 
		$\bigwedge_{(B_p,B_c)}  
		[
		{\{v^*(F_i): i\in B_p\} \not\subseteq\{1,\#_p\} 
		\vee
		\{v^*(F_i): i\in B_c\} \cap\{1,\#_c\}\not=\emptyset}
		\vee
		v^*(\Gamma)\not\subseteq\{1,\#_p\} 
		\vee
		v^*(\Delta)\cap\{1,\#_c\}\not=\emptyset
		]
		$

	iff 
		$\bigwedge_{(B_p,B_c)}
		[
		{\{v^*(F_i): i\in B_p\} \subseteq\{1,\#_p\} 
		\wedge
		v^*(\Gamma)\subseteq\{1,\#_p\} 
		]
		\Rightarrow
		[
		\{v^*(F_i): i\in B_c\} \cap\{1,\#_c\}\not=\emptyset}
		\vee
		v^*(\Delta)\cap\{1,\#_c\}\not=\emptyset
		]
		$	
}

\item[(iii)] From the two items above, using compactness and induction on the complexity of the (finite number of) formulae involved, we complete the proof that this semantics is sound and complete.
\qedhere
\end{itemize}
\end{myproof}

For regular and compact logics, Corollary~\ref{corol:maxSrank} about Suszko rank can now be subsumed by parallel results about the more stringent notion of truth-{functional} rank:

\begin{theorem}[Maximum truth-{functional} rank theorems]\label{thm:maxTCrank}
A monotonic, regular and compact logic is of truth-{functional} rank:
\begin{itemize}
	\item at most 2 if it is reflexive and transitive,
	\item at most 3 if it is transitive or reflexive,
	\item at most 4 in general.
\end{itemize}
\end{theorem}

\begin{myproof}
\begin{description}

\item[with reflexivity] 
To a monotonic reflexive logic, we apply Theorem~\ref{thm:monrefNOSI}, obtain a truth-adequate $p$-mixed semantics for it, and apply the {truth-functional} Scott-Suszko reduction.
We can be sure that no formula is assigned in any valuation $v$ the value $\#_p$. Otherwise, such a formula $F_p$ would violate reflexivity: $F_p\not\vdash F_p$ because in the relevant valuation $v_p$: $v_p(F_p)\subseteq\{1,\#_p\}$ and $v_p(F_p)\cap\{1,\#_c\}=\emptyset$. 
Hence, the {truth-functional} Scott-Suszko reduction is adequate, truth-interpretable and at most trivalent, over (some subset of) the set of truth values $\{1,\#_c,0\}$, which shows that monotonic, reflexive logics are of truth-{functional} rank at most 3.

\item[with transitivity]  If the logic is transitive, we can associate to it a truth-adequate $q$-mixed semantics by Theorem~\ref{thm:montransNOSI}, and take the {truth-functional} Scott-Suszko reduction of it. This semantics does not make use of $\#_c$, neither for the atoms, nor for more complex formulae: no $\#_c$ shows up in the atoms, that is in what we called $\mathcal{V}^*_v$ above, because $\mathcal{D}^\lambda_c\setminus\mathcal{D}^\lambda_p=\emptyset$; furthermore, if a $\#_c$ is generated for some formula $F_c$ in some valuation $v_c$, then that valuation can safely be dropped. If it could not be dropped, there would be $\Gamma\not\vdash\Delta$ such that $\forall v\not=v_c, v(\Gamma)\truthrelation^* v(\Delta)$. But then $\Gamma \not\vdash \Delta$ while both $\Gamma, F_c\vdash \Delta$ and $\Gamma \vdash F_c, \Delta$ would hold. This would violate transitivity.
Hence, we can obtain a truth-adequate mixed semantics over (a subset of) $\{1,\#_p,0\}$: the truth-compositional rank is at most 3.

\item[with reflexivity and transitivity] The previous construction for a transitive logic yields reflexivity if and only if $\#_p$ does not appear in any formula (a formula involving $\#_p$ would violate reflexivity, as explained above). The truth-{functional} rank is therefore at most 2. 

\item[without reflexivity nor transitivity] We apply the {truth-functional} Scott-Suszko reduction to a semantics given by Theorems~\ref{thm:monNOSI} to obtain a truth-adequate, mixed semantics over (a subset of) $\{1,\#_p,\#_c,0\}$. 
\qedhere
\end{description}
\end{myproof}

As before in Corollary~\ref{corol:exactSrank}, we can also obtain exact rank theorems for non-permeable logics:

\begin{theorem}[Exact truth-{functional} rank theorems]\label{thm:exactTCrank}
A non-permeable, monotonic, regular and compact logic is of mixed truth-{functional} rank:
\begin{itemize}
	\item exactly 2 if and only if it is reflexive and transitive,
	\item exactly 3 if and only if it is transitive or reflexive but not both.
	\item exactly 4 if and only if it is neither transitive nor reflexive.
\end{itemize}
\end{theorem}

\begin{myproof}
The proof is exactly the same as the proof for Corollary~\ref{corol:exactSrank}, except that the reference for the upper bounds is not Corollary~\ref{corol:maxSrank} but Theorem~\ref{thm:maxTCrank}, of course.
\end{myproof}

Hence, we obtained complete analogs of Corollaries~\ref{corol:maxSrank}/\ref{corol:exactSrank} for regular and compact logics in Theorems~\ref{thm:maxTCrank}/\ref{thm:exactTCrank}. We find that the Suszko rank is essentially the same as the \emph{a priori} more stringent truth-{functional} rank. Under these circumstances, then, logical and algebraic truth values are fundamentally the same, given appropriate reductions.\footnote{{As an example of nontrivial truth-functional reduction, we may give the case of \cite{smith:vdt}'s fuzzy logic. Smith proposes an infinite-valued semantics for the language of first-order logic and gives a definition of logical consequence which he shows to coincide with classical consequence. His proof is an instance of a Scott-Suszko reduction, which preserves truth-functionality. See \cite{cobreros2018tolerance} for discussion and further comparisons with Smith's logic.}} These results obtain only for compact and regular logics. If these restrictions are unsatisfying to some readers, Appendix~\ref{app:groupingreductions} explores a way to study more general logics, through what we call \emph{grouping reductions}, but we exhibit limitations of this other approach.

\section{Conclusions}

Let us summarize our progression in this paper. We have defined a \emph{logic} syntactically, as consisting of a set of formulae, a set of connectives, and a consequence relation definable as a set of arguments over that language. In parallel, we have defined a \emph{semantics} as a structure allowing us to map formulae, connectives, and the consequence relation to their interpretations. The rest of the paper has dealt with representation theorems, relating structural properties of a logic with general constraints on their semantic interpretation. 

Firstly, 
we have addressed Suszko's problem in its general form: what is the least number of truth values needed to characterize a logic semantically? We showed that if one is only interested in finding a sound and complete semantics, the least number of truth values needed to characterize a monotonic logic lies between 2 and 4, and the precise number is tied to whether the consequence relation is reflexive, transitive or neither of those. 
Similar results were {first} obtained by French and Ripley, using different techniques, {and by Blasio, Marcos and Wansing}. As the former emphasize, only partial correspondence results were obtained earlier (in particular by \citealp{humberstone1988heterogeneous}, \citealp{malinowski:q}, \citealp{tsuji1998many} and \citealp{frankowski2004}).

Secondly, we have pointed out that the rank theorems obtained no longer hold in full generality if we require the semantics of the connectives to be truth-{functional}. However, truth-{functionality} is a natural desideratum, and Suszko's problem can be extended as follows: what is the least number of truth values needed to represent a logic by means of a \emph{truth-{functional}} semantics? About this, we showed that for a natural class of logics, this generalization is not more demanding than Suszko's original question: the least number of truth values needed to truth-{functionally} characterize a monotonic, compact and regular logic is, again, exactly 2 if the logic is reflexive and transitive, exactly 3 if it is only one of those, exactly 4 if it is neither reflexive nor transitive.
This result is of significance for the understanding of what ought to count as a truth value.
Suszko pointed out that, given a logic, the number of algebraic truth values for a {truth-functional} semantics is in general arbitrary. On the other hand, he defended a notion of logical truth values, whose number is not arbitrary, but at the risk of losing {truth-functionality}. Our approach shows that we can have a notion of truth value that is tailored for truth-{functionality} and that may coincide with the logical notion.

Nevertheless, this result depends on the restriction to a special class of connectives and logics. To gain more generality, 
in Appendix~\ref{app:groupingreductions} we investigate  
whether more radical departures from the notion of reduction employed by Scott and Suszko could preserve {truth-functionality} including for noncompact nonregular logics. We introduce one possible candidate, the notion of grouping reduction. The notion is plausible from a semantic point of view in that it proposes to assimilate truth values that may be intersubstituted \emph{salva consequentia}. As we show there, however, it also leaves us with problems calling for further work. Because both the Scott-Suszko types of reduction and grouping reductions have limitations, more candidates may need to be considered to define an adequate notion of reduction, not only preserving consequence and {truth-functionality}, but ensuring that the reduction always helps reach an optimal number of truth values and may therefore reveal the rank of any logic and, in fact, a canonically associated semantic.

\bigskip

\appendix

\titleformat{\section}{\normalfont\Large\bfseries}{\appendixname~\thesection:}{1em}{}

\section*{Appendices}

\section{Simple compositionality in a Scott-Suszko reduction}\label{app:SSredcompos}

Compositionality, without truth-functionality, is a very weak demand. In particular, compositionality is met as soon as no two formulae are interpreted by the same proposition: there is then no real constraint on the interpretation of connectives at the level of propositions. This suggests two ways to artificially achieve compositionality. First, one may multiply the number of truth-values, even if they play roughly the same inferential role, just so that the truth values of formulae can differ in at least one world, even if in fact the truth values they are assigned are essentially the same. If one is looking for semantics with a minimal number of truth values, this is not the way to go, however. Another option is to distinguish semantically between formulae by multiplying the number of worlds: idle worlds may be added so that different formulae are assigned different truth values in these worlds, even if these worlds play no real inferential role. We pursue this second option in this Appendix.

To do so, we present a slightly different method to arrive at what is essentially the output of the original Scott-Suszko reduction. This construction is to be compared directly to the approach in \cite{french2017valuations}.
Interestingly, this construction is not a reduction \emph{per se}: it does in one step what was done in several steps in the core of the text (and in earlier discussions on the topic), that is, it exhibits a semantics with few truth values, not by reducing an existing semantics (typically obtained by the Lindenbaum method), but by extracting this semantics directly from the syntactic consequence relation. More to the point of this appendix, this method delivers a semantics in which each formula will most commonly be assigned to a unique proposition, because the semantics will use many worlds, and hence will be compositional (albeit not necessarily truth-{functional}).

\begin{definition}
The \emph{direct Scott-Suszko} semantics of a logic is defined as follows:

	$\mathcal{V}:=\{1,\#_p,\#_c,0\}$, 
	$\quad$
	$\mathcal{W}:=\{w_{(\Gamma, \Delta)} :\Gamma\not\vdash\Delta\}$,
	$\quad$	
	$\truthrelation:=\truthrelation_{\{1,\#_p\},\{1,\#_c\}}$,
	
	$\interp{F}(w_{(\Gamma, \Delta)}):=$
	{\scriptsize$\left\{\begin{array}{l@{\textrm{ if }}c@{\textrm{ and }}c}
	1 &  F\in\Gamma & F\not\in\Delta\\
	0 &  F\not\in\Gamma & F\in\Delta\\
	\#_p & F\in\Gamma & F\in\Delta\\
	\#_c &  F\not\in\Gamma & F\not\in\Delta
	\end{array}\right.$}

\end{definition}

\begin{theorem}
The direct Scott-Suszko semantics of a monotonic logic is sound and complete.
\end{theorem}

\begin{myproof}
If $\Gamma\not\vdash\Delta$, then it is easily shown that $\interp{\Gamma}(w_{(\Gamma,\Delta)})\subseteq\{1,\#_p\}$ and $\interp{\Delta}(w_{(\Gamma,\Delta)})\subseteq\{0,\#_p\}$. 
Conversely, if there is some world $w_{(\Gamma',\Delta')}$ such that 
$\interp{\Gamma}(w_{(\Gamma',\Delta')})\not\truthrelation\interp{\Delta}(w_{(\Gamma',\Delta')})$, then
$\Gamma\subseteq\Gamma'$ and $\Delta\subseteq\Delta'$. Hence, by monotonicity, $\Gamma\not\vdash\Delta$.
\end{myproof}

\begin{theorem}
If there is at most one formula $F$ such that $\forall\Gamma,\Delta:\Gamma,F\vdash\Delta \textrm{ and } \Gamma\vdash F, \Delta$, then the direct Scott-Suszko semantics of a monotonic logic assigns a different proposition to all formulae.
\end{theorem}

\begin{myproof}
Let $F$ be a formula such that it is not the case that $\forall\Gamma,\Delta:\Gamma,F\vdash\Delta \textrm{ and } \Gamma\vdash F, \Delta$. By monotonicity, this amounts to $F\not\vdash$ or $\not\vdash F$. If the former is the case, then $F$ is the only formula such that $\interp{F}(w_{(\{F\},\emptyset)})\subseteq\{1,\#_p\}$; if the latter is the case, then $F$ is the only formula such that $\interp{F}(w_{(\emptyset,\{F\})})\subseteq\{0,\#_p\}$.
\end{myproof}

Compositionality is a very weak demand, it holds in particular if all formulae are assigned a different proposition. Thus, compositionality of the direct Scott-Suszko semantics is an immediate corollary:

\begin{corollary}
If every formula $F$ is involved in at least one relation that does not hold ($\exists\Gamma,\Delta: F\in(\Gamma\cup\Delta) \textrm{ and } \Gamma\not\vdash\Delta$),
the direct Scott-Suszko semantics of a monotonic logic is compositional.
\end{corollary}

The reader can see why from this construction we obtain (mixed) compositional rank results similar to the (mixed) Suszko rank results: certainly compositional rank is at most 4, and reflexivity would lead to drop $\#_p$ and transitivity can help drop $\#_c$ as a truth value (although it is a little more complicated for transitivity, which requires to take a subset of $W$ and therefore lose some of the leverage to assign a unique proposition to each formula).

\section{The canonical truth-functions of classical connectives}\label{app:regconnectives}

The {truth-functional} version of the Scott-Suzsko reduction (Definition~\ref{def:strongssreduction}) introduces a canonical way to associate truth-functions to regular connectives, given their regularity rules. For usual connectives from, say, bivalent logic, which come with regularity rules, these truth-functions offer a canonical 4-valued extension of their bivalent truth-function. Below we show the truth-tables associated in this way to the connectives mentioned in Example~\ref{def:structconn}, see the regularity rules reported up there. More connectives are shown with their regularity rules and canonical truth-tables in the semi-automatic Excel file attached with this paper as an ancillary document.\medskip

\begin{tabular}{cccc}
Negation &
Conjunction &
Disjunction &
Conditional\\

\framebox{\begin{tabular}{c|c}
	 & $\neg$\\ \hline
$1$	&$0$\\
$\#_p$	&$\#_p$\\
$\#_c$	&$\#_c$\\
$0$	&$1$\end{tabular}}

& 
\framebox{\begin{tabular}{c|cccc}
	 & $1$ 		& $\#_p$ 	& $\#_c$ 	&$0$\\ \hline
$1$	&$1$		&$\#_p$	&$\#_c$	&$0$\\
$\#_p$	&$\#_p$	&$\#_p$	&$0$		&$0$\\
$\#_c$	&$\#_c$	&$0$		&$\#_c$	&$0$\\
$0$	&$0$	&$0$	&$0$		&$0$\end{tabular}}

& 
\framebox{\begin{tabular}{c|cccc}
	 &$1$ 	& $\#_p$ 	& $\#_c$ 	&$0$\\ \hline
$1$	&$1$	&$1$		&$1$		&$1$\\
$\#_p$	&$1$	&$\#_p$	&$1$		&$\#_p$\\
$\#_c$	&$1$	& $1$		&$\#_c$	&$\#_c$\\
$0$	&$1$	&$\#_p$	&$\#_c$		&$0$\end{tabular}}

& \framebox{\begin{tabular}{c|cccc}
	 &$1$ 	& $\#_p$ 	& $\#_c$ 	&$0$\\ \hline
$1$	&$1$	&$\#_p$	&$\#_c$	&$0$\\
$\#_p$	&$1$	&$\#_p$	&$1$		&$\#_p$\\
$\#_c$	&$1$	&$1$		&$\#_c$	&$\#_c$\\
$0$	&$1$	&$1$		&$1$		&$1$\end{tabular}}

\end{tabular}\medskip

How do the above truth-tables compare with other truth-tables proposed for many-valued negation, conjunction, conditionals, ...?
Below we highlight two properties of such truth-tables.

\begin{definition}[Closure]\label{def:stability}
A set of truth values $\underline{\mathcal{V}}\subseteq\mathcal{V}$ is \emph{closed} under a truth-function $\sem{C}$ iff
\[\forall \alpha_1, ..., \alpha_n\in\underline{\mathcal{V}}: C(\alpha_1, ..., \alpha_n)\in\underline{\mathcal{V}}\]
\end{definition}

\begin{definition}\label{def:bivalentconnective}
A regular connective $C$ is called \emph{bivalent-closed} if $\{0,1\}$ is closed under its canonically associated truth-function $\sem{C}$.
\end{definition}

\begin{theorem}\label{thm:stability}
For a bivalent-closed connective $C$, $\{0,1,\#_p\}$ and $\{0,1,\#_c\}$ are closed under $\sem{C}$.
\end{theorem}

\begin{myproof}
Assume that $\{0,1\}$ is closed and that $\alpha_1, ..., \alpha_n\subseteq\{0,1,\#_p\}$. Assume that $C(\alpha_1, ..., \alpha_n)=\#_c$, then 

	$\sem{C}(\alpha_1, ..., \alpha_n)\not\in\{1,\#_p\}$ so
		$\bigwedge\limits_{(B_p,B_c)\in \mathcal{B}^p} 
			{(\{\alpha_i: i\in B_p\} \subseteq\{1,\#_p\} 
			\Rightarrow 
			\{\alpha_i: i\in B_c\} \cap\{1,\#_c\}\not=\emptyset)}$
	
	$\sem{C}(\alpha_1, ..., \alpha_n)\in\{1,\#_c\}$ so
		$\bigwedge\limits_{(B_p,B_c)\in \mathcal{B}^c} 
			{(\{\alpha_i: i\in B_p\} \subseteq\{1,\#_p\} 
			\Rightarrow 
			\{\alpha_i: i\in B_c\} \cap\{1,\#_c\}\not=\emptyset)}$

\noindent
None of these two statements would become false if some $\alpha_i$ was changed from $\#_p$ to $1$. Hence, we would have $C(\pi(\alpha_1), ..., \pi(\alpha_n))=\#_c$, with $\pi(0):=0$, $\pi(1):=1$ and $\pi(\#_p):=1$. And this would violate the stability of $\{0,1\}$. This proves the closure of $\{1,0,\#_p\}$. 
The proof for the closure of $\{1,0,\#_c\}$ is similar: no conjunct in the statement above can become true if some $\alpha_i$ is changed from $\#_c$ to $0$.
\end{myproof}

As a result, if we ignore one of the non-classical truth value as an argument, $\#_p$ or $\#_c$, we are left with a trivalent connective in its output too. In other words, removing the $\#_p$ row and the $\#_p$ column in the truth-table above for instance, leaves us with a truth-table which does not contain $\#_p$ anymore. On these three-valued spaces then, either $\{1,\#_c,0\}$ or $\{1,\#_p,0\}$, we recognize the usual trivalent Strong Kleene tables, in the following sense:

\begin{definition}[Strong Kleene]\label{def:strongkleene}
A truth-function is \emph{Strong Kleene} if for $\alpha_1, ..., \alpha_{k-1}, \alpha_{k+1}, ..., \alpha_{n}\subseteq\{0,1\}$,
	whenever $\sem{C}(\alpha_1, ..., \alpha_{k-1}, 0, \alpha_{k+1}, ..., \alpha_{n})$
	and
	$\sem{C}(\alpha_1, ..., \alpha_{k-1}, 1, \alpha_{k+1}, ..., \alpha_{n})$
	take the same value, then
	$\sem{C}(\alpha_1, ..., \alpha_{k-1}, \alpha, \alpha_{k+1}, ..., \alpha_{n})$ takes this same value for all truth values $\alpha$.
\end{definition}

Interestingly, the canonical truth-functions of the 16 classical bivalent binary connectives are Strong Kleene, as can be seen from the computations of all these truth-functions in the auxiliary file:

\begin{fact}
The truth-functions canonically associated to the 16 classical binary connectives {are} Strong Kleene.
\end{fact}

It would be useful to investigate a structural characterization of bivalent-stable connectives, of Strong Kleene connectives and, in connection to this, to ask when several regularity rules may co-exist for a connective (it certainly is so for `xor') because they could lead to different `canonical' truth-functions (although there is a canonical way to decide between several regularity rules, by picking the conjunction of them all, which remains a regularity rule). We leave such a systematic investigation of connectives for another occasion.

\section{Grouping reductions}\label{app:groupingreductions}

{This section presents a notion of reduction which is natural from an algebraic perspective, but in fact less efficient than the Suszko reduction, for reasons explained hereafter.}\footnote{{Some of the results below may be derived in a more abstract setting. We are indebted to a referee for pointing out that Theorem \ref{thm:redequivclass} in particular is an instance of a more general result on Leibniz congruences, see in particular \cite{blok1989algebraizable}, \cite{font1991leibniz}.}}

\subsection{Definition of grouping reductions}

Can we study more stringent notions of rank for non-compact and non-regular logics? To do so, one may like to exhibit a reduction that preserves truth-adequacy, and reaches the truth-compositional rank --- just like Scott-Suszko reductions preserve truth-relationality and reach their associated rank under appropriate circumstances (see Corollary~\ref{corol:exactSrank} and Theorem~\ref{thm:exactTCrank}).
With this in hand, we would have a canonical way to construct an optimized (algebraic) semantics for any logic: start from any semantics that can be obtained for a logic, e.g., the Lindenbaum construction we described, and apply the reduction. 

A natural way to reduce the number of truth values is to group {functionally} similar truth values together. In other words, one could form new truth values as equivalence classes of old truth values. Formally, we may project the set of truth values into a smaller set of truth values, essentially collapsing truth values with the same image in this new set. We propose a formal description of this process, exhibiting constraints to preserve the soundness and completeness of the relation on the one hand, and compositionality for connectives on the other hand.

\begin{definition}[Grouping reductions and associated properties]
Consider a logic and an intersective semantics for it over the set of truth values $\mathcal{V}_1$, with truth-relation $\truthrelation_1$, and connectives generically called $C$ interpreted as $\sem{C}_1$. 
A surjection $\rho:\mathcal{V}_1\twoheadrightarrow\mathcal{V}_2$ is called a \emph{grouping-reduction} (abbreviated \emph{g-reduction}) and can be qualified:

\begin{itemize}

\item $\rho$ is a \emph{relation-g-reduction} if
$\forall\gamma, \delta: \forall \gamma', \delta':$ \\
	$$\textrm{if }(\rho(\gamma)=\rho(\gamma') \textrm{ and } \rho(\delta)=\rho(\delta')),
	\quad\textrm{ then }
	(\gamma\truthrelation_1\delta \textrm{ iff } \gamma'\truthrelation_1\delta')
	$$ 

\item $\rho$ is a \emph{$C$-g-reduction} for some connective $C$ if
$\forall x_1, ..., x_n: \forall x'_1, ..., x'_n$: \\
	$$\textrm{if }(\rho(x_1)=\rho(x'_1), ..., \rho(x_n)=\rho(x'_n)), 
	\quad\textrm{ then }
	(\rho(\sem{C}_1(x_1, ..., x_n))=\rho(\sem{C}_1(x'_1, ..., x'_n)))$$ 

\item $\rho$ is a \emph{strong g-reduction} if it is both a relation-g-reduction and a $C$-g-reduction for all connectives.

\end{itemize}

\end{definition}

The above definition articulates a constraint concerning the truth-relation and a constraint for connectives. The former is in line with Scott and Suszko's approach (although the constraint here is more stringent, as we will see). The latter however builds into the enterprise the hope that compositionality will be maintained. The following theorem, indeed, shows how a strong g-reduction provides a way to define a truth-relational semantics with $|\mathcal{V}_2|$-many truth values. We will also call a semantics defined in this way a \emph{strong g-reduction}:
	
\begin{theorem}
Consider a logic, with a strong intersective mixed semantics for it (over the set of truth values $\mathcal{V}_1${)} and a strong g-reduction of it $\rho:\mathcal{V}_1\twoheadrightarrow\mathcal{V}_2$.
A strong, intersective mixed semantics for the logic can be defined over $\mathcal{V}_2$.
\end{theorem}

\begin{myproof}
We exhibit a valuational, compositional, truth-relational and truth-functional semantics over $\mathcal{V}_2$ by defining its interpretation function inductively (see Fact~\ref{fact:semanticsfromtruth}). Pick an inverse function $\rho^{\langle-1\rangle}$ for the surjective function $\rho$ (i.e.~$\rho\circ\rho^{\langle-1\rangle}=\textrm{Id}$). 
Then define 
$\truthrelation_2:=\bigcap\truthrelation_{\rho(\mathcal{D}_p^{\lambda}),\ \rho(\mathcal{D}_c^{\lambda})}$
and $\sem{C}_2:=\rho\circ\sem{C}_1\circ\rho^{\langle-1\rangle}$. 

The core of the proof is to show that this semantics is adequate. This holds because of the equivalence between 1 and 8 below. In fact, we will prove that all the following statements are equivalent: 

\begin{enumerate}
\item $\Gamma\vdash\Delta$
\item $\forall v_1: (\mathcal{A}\to\mathcal{V}_1)\quad
	 v_1 (\Gamma)\truthrelation_1 v_1 (\Delta)$
\item $\forall v_1: (\mathcal{A}\to\mathcal{V}_1)\quad
	\forall \gamma, \delta [\rho(\gamma)= \rho(v_1 (\Gamma)) \textrm{ and } \rho(\delta)= \rho(v_1 (\Delta))]:
		\gamma\truthrelation_1\delta$
\item $\forall v_1: (\mathcal{A}\to\mathcal{V}_1)\quad
	\rho^{-1}(\rho( v_1 (\Gamma)))\truthrelation_1\rho^{-1}(\rho( v_1 (\Delta)))$
\item $\forall v_1: (\mathcal{A}\to\mathcal{V}_1)\quad
	\forall\lambda: 
	\rho^{-1}(\rho( v_1 (\Gamma)))\subseteq\mathcal{D}_p^{\lambda}
	\Rightarrow
	\rho^{-1}(\rho( v_1 (\Delta)))\cap\mathcal{D}_c^{\lambda}\not=\emptyset$
\item $\forall v_1: (\mathcal{A}\to\mathcal{V}_1)\quad
	\forall \gamma, \delta [\rho(\gamma)= v_1 (\Gamma) \textrm{ and } \rho(\delta)= v_1 (\Delta)]:
		\gamma\subseteq\mathcal{D}_p^{\lambda}
		\Rightarrow
		\delta\cap\gamma\subseteq\mathcal{D}_c^{\lambda}\not=\emptyset$

\item $\forall v_1: (\mathcal{A}\to\mathcal{V}_1)\quad
	\forall\lambda: 
	\rho( v_1 (\Gamma))\subseteq\rho(\mathcal{D}_p^{\lambda})
	\Rightarrow
	\rho( v_1 (\Delta))\cap\rho(\mathcal{D}_c^{\lambda})\not=\emptyset$

\item $\forall v_2: (\mathcal{A}\to\mathcal{V}_2)\quad
	\forall\lambda: 
	 v_2 (\Gamma)\subseteq\rho(\mathcal{D}_p^{\lambda})
	\Rightarrow
	 v_2 (\Delta)\cap\rho(\mathcal{D}_c^{\lambda})\not=\emptyset$

\item $\forall v_2: (\mathcal{A}\to\mathcal{V}_2)\quad
	 v_2 (\Gamma)\truthrelation_2 v_2 (\Delta)$

\end{enumerate}

\noindent
\textbf{1$\Leftrightarrow$2} is the adequacy of semantics 1.
\textbf{2$\Leftrightarrow$3} is because $\rho$ is a relation-reduction.
\textbf{2$\Leftrightarrow$4} is because $\rho$ is a relation-reduction and $\rho(\rho^{-1}(X))=X$ for all $X$ since $\rho$ surjective.
\textbf{4$\Leftrightarrow$5} is the definition of $\truthrelation_1$ as an intersection of mixed consequence relations.
\textbf{3$\Leftrightarrow$6} is the definition of $\truthrelation_1$ as an intersection of mixed consequence relations.
\textbf{8$\Leftrightarrow$9} is the definition of $\truthrelation_2$ as an intersection of mixed consequence relations.
\begin{description}
\item[7$\Leftrightarrow$8] is because $ v_2 (F)=\rho( v_1 (F))$, for all $v_1$ such that $\forall p: v_2(p)=\rho(v_1(p))$ (note this $v_1\sim_{\rho/\mathcal{A}} v_2$). This is shown by induction. It is clear if $F$ is an atom. And then:\\
\begin{tabular}{l@{~=~}lll}
$v_2(C(F_1, ..., F_n)$ 
	& $\rho \circ \sem{C}_1 \circ \rho^{\langle-1\rangle} (v_2(F_1), ..., v_2(F_n))$\\
	& $\rho \circ \sem{C}_1 \circ \rho^{\langle-1\rangle} [ \rho (v_1(F_1) ), ..., \rho (v_1(F_n) ) ]$
		& for all $v_1\sim_{\rho/\mathcal{A}} v_2$
		& (by induction)\\
	& $\rho \circ \sem{C}_1 [ v_1(F_1), ..., v_1(F_n) ]$
		& for all $v_1\sim_{\rho/\mathcal{A}} v_2$
		& (by $C$-reduction)\\
	& $\rho ( v_1 (C(F_1, ..., F_n)) )$
		& for all $v_1\sim_{\rho/\mathcal{A}} v_2$
		& (by truth-functionality)\\
\end{tabular}

\item[6$\Rightarrow$7] 
	Suppose 7 does not hold, then pick some $v_1$, $\lambda$ such that 
	$\rho( v_1 (\Gamma))\subseteq\rho(\mathcal{D}_p^{\lambda})$
	and
	$\rho( v_1 (\Delta))\cap\rho(\mathcal{D}_c^{\lambda})=\emptyset$.
	Then define 
		$\gamma:=\rho^{-1}(\rho( v_1 (\Gamma)))\cap\mathcal{D}_p^{\lambda}$ 
		and 
		$\delta:=\rho^{-1}(\rho( v_1 (\Delta)))\setminus\mathcal{D}_c^{\lambda}$. 
	We obtain a violation of 6 because:
	\begin{description}
	\item[$\gamma\subseteq\mathcal{D}_p^{\lambda}$ and $\delta\cap\mathcal{D}_c^{\lambda}=\emptyset$:] by construction.
	\item[$\rho(\gamma)=\rho( v_1 (\Gamma))$:]
		$\rho(\gamma) 
			= \{\rho(x): x\in\gamma\}
			= \{\rho(x): x\in\rho^{-1}(\rho( v_1 (\Gamma))) \textrm{ and }x\in\mathcal{D}_p^{\lambda}\}
			= \{\rho(x): \rho(x)\in\rho( v_1 (\Gamma)) \textrm{ and }x\in\mathcal{D}_p^{\lambda}\}
			= \rho( v_1 (\Gamma))\cap\rho(\mathcal{D}_p^{\lambda})
			= \rho( v_1 (\Gamma))
		$
	\item[$\rho(\delta)=\rho( v_1 (\Delta))$:]
		$\rho(\delta) 
			= \{\rho(x): x\in\delta\}
			= \{\rho(x): x\in\rho^{-1}(\rho( v_1 (\Delta))) \textrm{ and }x\not\in\mathcal{D}_c^{\lambda}\}
			= \{\rho(x): \rho(x)\in\rho( v_1 (\Delta)) \textrm{ and }x\in(\mathcal{V}_1\setminus\mathcal{D}_c^{\lambda})\}
			= \rho( v_1 (\Delta))\cap(\rho(\mathcal{V}_1\setminus\mathcal{D}_c^{\lambda}))
			= \rho( v_1 (\Delta))
		$
	\end{description}
\item[7$\Rightarrow$5] 
	Assume that 5 does not hold: for some $\lambda,  v_1 $: 
	(i)~$\rho^{-1}(\rho( v_1 (\Gamma)))\subseteq\mathcal{D}_p^{\lambda}$
	and
	(ii)~$\rho^{-1}(\rho( v_1 (\Delta)))\cap\mathcal{D}_c^{\lambda}=\emptyset$.
	Applying $\rho$ to both sides of (i), using $\rho(\rho^{-1}(X))=X$, we obtain $\rho( v_1 (\Gamma)))\subseteq\rho(\mathcal{D}_p^{\lambda})$.
	Assume that there is $y$ in $\rho( v_1 (\Delta))\cap\rho(\mathcal{D}_c^{\lambda})$. Then, there is $x\in\mathcal{D}_c^{\lambda}$ such that $\rho(x)\in\rho( v_1 (\Delta))$, that is $x\in\rho^{-1}(\rho( v_1 (\Delta))\cap\mathcal{D}_c^{\lambda}$. This would contradict (ii). 
\qedhere\end{description}

\end{myproof}

\subsection{Reduction and equivalence classes of truth values}\label{app:equivalencered}

G-reductions group together truth values that play the same role. The truth values after reduction are classes of truth values, grouped together based on their image by $\rho$. Here we show that under certain circumstances, the optimal equivalence relation can be inferred from the truth-relation.

\begin{definition}\label{def:redequivclass}
Given a semantics, we define an equivalence relation over truth values as 
	$$x\sim y
	\quad\textrm{ iff }\quad
	(\forall\gamma,\delta: \gamma, x\truthrelation_1\delta \textrm{ iff } \gamma, y\truthrelation_1\delta) 
	\textrm{ and }
	(\forall\gamma,\delta: \gamma\truthrelation_1 x, \delta \textrm{ iff } \gamma\truthrelation_1 y, \delta)
	$$
The surjection from truth values to their equivalence class is called the canonical g-reduction, noted~$\widetilde{\rho}$.
\end{definition}

The following theorem shows that relation-g-reductions can only collapse truth values that the canonical g-reduction collapses. The canonical g-reduction therefore imposes a lower bound on the number of truth values a relation-g-reduction may involve.

\begin{theorem}\label{thm:max-g-red}
For any relation-g-reduction $\rho$, if $\rho(x)=\rho(y)$ then $\widetilde{\rho}(x)=\widetilde{\rho}(y)$. 
\end{theorem}
	
\begin{myproof}
Suppose $\rho(x)=\rho(y)$. 
Then $\forall \gamma, \delta: 
	\rho(\gamma,x)=\rho(\gamma,y)$ 
	and 
	$\rho(\delta,x)=\rho(\delta,y)$.
Since $\rho$ is a relation-reduction, it follows that:
	($\gamma, x\truthrelation\delta$ iff $\gamma, y\truthrelation\delta$)
	and
	($\gamma \truthrelation x, \delta$ iff $\gamma \truthrelation y, \delta$),
i.e.~$x \sim y$,
i.e.~$\widetilde{\rho}(x)=\widetilde{\rho}(y)$.
\end{myproof}

But this lower bound is not necessarily reached, because the canonical g-reduction is not in general a relation-g-reduction (Example~\ref{ex:not-relation-red}). Theorem~\ref{thm:redequivclass}, however, shows that with some finiteness assumptions it is a relation-g-reduction (Theorem~\ref{thm:redequivclass}).

	\begin{example}\label{ex:not-relation-red}
	The canonical g-reduction is not in general a relation-g-reduction. 	\end{example}
	
	\begin{myproof}
	Consider $\mathcal{V} = \{0, 1\} \cup \{\#_n : n \in N\}$, and $\truthrelation=\bigcap_{\mathcal{D} \textrm{ finite subset of } \mathcal{V}\setminus\{0\}}\truthrelation_{\mathcal{D},\{1\}} $. For all $i,j$, $\#_i\sim\#_j$. But the following shows that this does not yield a reduction:
	{for} $\widetilde{\{\#_n : n \in N\}}=\widetilde{\{\#_1\}}$,
	{and} $\{\#_n : n \in N\}\truthrelation 0$, but
	$\#_1\not\truthrelation 0$.
	In other words, one cannot collapse all the $\#_n$s into finitely many values $\#$.
	\end{myproof}

\begin{theorem}\label{thm:redequivclass}
When a semantics has a relation-g-reduction with a finite number of truth values, the canonical g-reduction is a relation-g-reduction, one using the fewest possible truth values given Theorem~\ref{thm:max-g-red}.
\end{theorem}

\begin{myproof}
Consider $\gamma, \delta, \gamma', \delta'$ such that $\widetilde{\rho}(\gamma)=\widetilde{\rho}(\gamma')$, $\widetilde{\rho}(\delta)=\widetilde{\rho}(\delta')$. We will show that $\gamma\truthrelation_1\delta\Leftrightarrow\gamma'\truthrelation_1\delta'$ by assuming one of them arbitrarily, say $\gamma\truthrelation_1\delta$ and show that the other follows.
\begin{itemize}

\item Assume first that $\gamma, \delta, \gamma', \delta'$ are finite. By monotonicity, we note that $\gamma,\gamma'\truthrelation_1\delta, \delta'$. From there, we can replace elements of $\gamma$ and $\delta$ one by one with equivalent elements in $\gamma'$ and $\delta'$ respectively, until $\gamma\cup\gamma'$ is reduced to $\gamma'$ and $\delta\cup\delta'$ is reduced to $\delta'$. The equivalence relation guarantees that the relation is preserved at each step and since there are finitely many replacements, we obtain $\gamma'\truthrelation_1\delta'$.

\item Consider a g-reduction $\rho:\mathcal{V}_1\twoheadrightarrow\mathcal{V}_2$ with $\mathcal{V}_2$ finite, and an `inverse' function for the surjection: $\rho^{\langle-1\rangle}$ such that $\rho\circ\rho^{\langle-1\rangle}=\textrm{Id}$.
For $\gamma, \delta\subseteq\mathcal{V}_1$, consider
	$\overline{\gamma}:=\rho^{\langle-1\rangle}\circ\rho(\gamma)$ and
	$\overline{\delta}:=\rho^{\langle-1\rangle}\circ\rho(\delta)$. 
	We obtain that:
	\begin{itemize}
	\item $\overline{\gamma}$ and $\overline{\delta}$ are finite, since $\mathcal{V}_2$ is finite.
	\item ($\gamma\truthrelation_1\delta$ iff $\overline{\gamma}\truthrelation_1\overline{\delta}$), because $\gamma$ and $\delta$ and $\overline{\gamma}$ and $\overline{\delta}$ have the same images by $\rho$, which is a relation g-reduction.
	\item $\widetilde{\rho}(\overline{\gamma})=\widetilde{\rho}(\gamma)$ and $\widetilde{\rho}(\overline{\delta})=\widetilde{\rho}(\delta)$ (by Theorem~\ref{thm:max-g-red}). 
	\end{itemize}
	For all $\gamma, \gamma', \delta, \delta'$ we can thus reason with their $\overline{\gamma}, \overline{\gamma'}, \overline{\delta}, \overline{\delta'}$ finite counterparts and apply the argument from the previous step.
\qedhere
\end{itemize}

\end{myproof}

\subsection{Properties of grouping reductions}
	
Strong g-reductions are defined by constraining how truth values are grouped together depending on the role they play for the truth-relation (relation-g-reduction) and for each connective $C$ ($C$-g-reduction). We can provide an example showing how strong g-reductions differ from the original Scott-Suszko reduction. In Theorem~\ref{thm:notrured-n}, we showed that the Scott-Suszko reduction {can fail} to provide a truth-adequate semantics when we start from an \emph{order-theoretic semantics}, that is a semantics for which the truth-relation is essentially obtained by a total order on the set of truth values and the rule $\gamma\truthrelation\delta$ iff $\textrm{Inf}(\gamma)\leq\textrm{Sup}(\delta)$. A strong g-reduction would not fall into this problem, simply because there is no strict strong g-reduction of an order-theoretic semantics.
	
\begin{theorem}\label{thm:ordertheoreticnotreducible}
A strong order-theoretic semantics has no strong g-reduction.
\end{theorem}
	
\begin{myproof}
Consider a strong order-theoretic semantics, that is a strong semantics of which the truth-relation is derived from a well-order~$\leq$ on truth values as $\gamma\truthrelation\delta$ iff (def) $\textrm{Inf}(\gamma)\leq\textrm{Sup}(\delta)$. Then for any two distinct truth values, $x\truthrelation x$ but either $x\not\truthrelation y$ or $y\not\truthrelation x$. Hence, $x\not\sim y$, and Theorem~\ref{thm:max-g-red} guarantees that no g-reduction would collapse $x$ and $y$.
\end{myproof}

Hence, strong g-reductions do not reduce a semantics to the point of destroying compositionality. But g-reductions may be too weak: one may like to ensure that strong g-reductions reach the strong rank. This further desideratum is not satisfied:
	
\begin{theorem}\label{thm:notreachrank}
Some strong semantics may not be g-reducible to a semantics for which the cardinal of the set of truth values is the strong rank of the logic.
\end{theorem}
	
\begin{myproof}
Consider propositional logic and its classical semantics. Now consider the product semantics: $\mathcal{V}_*=\mathcal{V}\times\mathcal{V}$, $\interp{F}_*:=(\interp{F},\interp{F})$, $(\gamma_1,\gamma_2) \truthrelation_* (\delta_1, \delta_2)$ iff (def) ($\gamma_1\truthrelation_{\textrm{\textsc{biv}}} \delta_1$ and $\gamma_2\truthrelation_{\textrm{\textsc{biv}}} \delta_2$), with $\truthrelation_{\textrm{\textsc{biv}}}$ the classical bivalent truth-relation $\truthrelation_{\{1\},\{1\}}$. We show that (i)~this semantics is strong (the main issue is whether it is adequate), (ii)~this semantics cannot be g-reduced.
\begin{enumerate}
\item[(i)] 
	$\forall v_*: v_*(\Gamma)\truthrelation_* v_*(\Delta)$ 

	iff $\forall v_1, v_2: (v_1(\Gamma), v_2(\Gamma))\truthrelation_* (v_1(\Delta), v_2(\Delta))$ 

	iff $\forall v_1, v_2: (v_1(\Gamma) \truthrelation_1 v_1(\Delta)$ and $v_2(\Gamma) \truthrelation_2 v_2(\Delta))$ 

	iff ($\forall v_1: v_1(\Gamma)\truthrelation_1 v_1(\Delta)$ 
		and $\forall v_2: v_2(\Gamma)\truthrelation_2 v_2(\Delta)$)

	iff $\Gamma\vdash\Delta$ 

\item[(ii)] None of the values $(0,0)$, $(0,1)$, $(1,0)$, $(1,1)$ are equivalent (see Theorem~\ref{thm:redequivclass}), hence there is no strict relation-g-reduction of this semantics.
\end{enumerate}
It is important to note that over $\mathcal{V}\times\mathcal{V}$, one could also define the truth-relation as $\gamma\truthrelation\delta$ iff (def) $\pi_1(\gamma)\truthrelation\pi_2(\delta)$, ignoring the second member of each pair, and yet obtain an adequate semantics. The semantics with this truth-relation could be g-reduced, however, even though it is essentially the same. The theorem obtains because \emph{some} semantics cannot be g-reduced, but it suggests an extension of our notion of reduction, starting by a reduction of the representation as intersection of mixed consequence relations, but first dropping any member of the intersection that can be dropped, and then applying a g-reduction.
\qedhere
\end{myproof}

Hence, the constraint is too stringent in that g-reductions do not always allow us to reach the rank (Theorem~\ref{thm:notreachrank}). 
If we give up the Suszko rank in favor of the truth-compositional rank or strong rank, and want to study not only compact and regular logics but all logics, future research is needed to find a constructive reduction that can reveal the minimal number of truth values for a truth-relational \emph{and} (truth-)compositional semantics. This appendix shows that it will involve a deeper re-organization of truth values than simply grouping them together.

\bibliographystyle{apalike}
\bibliography{chemla_egre_rsl}

\newcommand{\SortNoop}[1]{}
\begin{thebibliography}{}

\bibitem[Andr{\'e}ka et~al., 2001]{andreka2001algebraic}
Andr{\'e}ka, H., N{\'e}meti, I., and Sain, I. (2001).
\newblock Algebraic logic.
\newblock In Gabbay, D. and Guenthner, F., editors, {\em Handbook of
  philosophical logic}, volume~2, pages 133--247. Springer.

\bibitem[B{\'e}ziau, 2001]{beziau2001sequents}
B{\'e}ziau, J.-Y. (2001).
\newblock Sequents and bivaluations.
\newblock {\em Logique et Analyse}, 44(176):373--394.

\bibitem[Blasio et~al., 2018]{blasio2017inferentially}
Blasio, C., Marcos, J., and Wansing, H. (2018).
\newblock An inferentially many-valued two-dimensional notion of entailment.
\newblock {\em Bulletin of the Section of Logic}, 46:233--262.

\bibitem[Blok and Pigozzi, 1989]{blok1989algebraizable}
Blok, W.~J. and Pigozzi, D. (1989).
\newblock {\em Algebraizable logics}, volume~77.
\newblock American Mathematical Society.

\bibitem[Bloom et~al., 1970]{bloom1970theorems}
Bloom, S.~L., Brown, D.~J., and Suszko, R. (1970).
\newblock Some theorems on abstract logics.
\newblock {\em Algebra and Logic}, 9(3):165--168.

\bibitem[Caleiro et~al., 2003]{caleiro2003dyadic}
Caleiro, C., Carnielli, W.~A., Coniglio, M.~E., and Marcos, J. (2003).
\newblock Dyadic semantics for many-valued logics.
\newblock Preprint available at: \url{http://wslc. math. ist. utl.
  pt/ftp/pub/CaleiroC/03-CCCM-dyadic2. pdf}.

\bibitem[Caleiro and Marcos, 2012]{caleiro2012many}
Caleiro, C. and Marcos, J. (2012).
\newblock Many-valuedness meets bivalence: Using logical values in an effective
  way.
\newblock {\em Multiple-Valued Logic and Soft Computing}, 19(1-3):51--70.

\bibitem[Caleiro et~al., 2015]{caleiro2015bivalent}
Caleiro, C., Marcos, J., and Volpe, M. (2015).
\newblock Bivalent semantics, generalized compositionality and analytic
  classic-like tableaux for finite-valued logics.
\newblock {\em Theoretical Computer Science}, 603:84--110.

\bibitem[Chemla and Egr\'e, 2019]{chemla2019mv}
Chemla, E. and Egr\'e, P. (2019).
\newblock From many-valued consequence to many-valued connectives.
\newblock Manuscript.

\bibitem[Chemla et~al., 2017]{CES2016}
Chemla, E., Egr\'e, P., and Spector, B. (2017).
\newblock Characterizing logical consequence in many-valued logic.
\newblock {\em Journal of Logic and Computation}, 27(7):2193--2226.

\bibitem[Cobreros et~al., 2012a]{cobreros:mixed}
Cobreros, P., Egr{\'e}, P., Ripley, D., and van Rooij, R. (2012a).
\newblock Tolerance and mixed consequence in the s'valuationist setting.
\newblock {\em Studia logica}, 100(4):855--877.

\bibitem[Cobreros et~al., 2012b]{tcs}
Cobreros, P., Egr\'e, P., Ripley, D., and van Rooij, R. (2012b).
\newblock Tolerant, classical, strict.
\newblock {\em Journal of Philosophical Logic}, 41(2):347--385.

\bibitem[Cobreros et~al., 2018]{cobreros2018tolerance}
Cobreros, P., Egr\'e, P., Ripley, D., and van Rooij, R. (2018).
\newblock Tolerance and degrees of truth.
\newblock Manuscript.

\bibitem[Font, 1991]{font1991leibniz}
Font, J.~M. (1991).
\newblock On the {L}eibniz congruences.
\newblock In {\em Algebraic Methods in Logic and Computer Science}, volume~28,
  pages 17--34. Banach Center Publications.

\bibitem[Font, 2003]{font2003generalized}
Font, J.~M. (2003).
\newblock Generalized matrices in abstract algebraic logic.
\newblock In Hendricks, V. and Malinowski, J., editors, {\em Trends in logic},
  volume~21, pages 57--86. Springer.

\bibitem[Font, 2009]{font2009taking}
Font, J.~M. (2009).
\newblock Taking degrees of truth seriously.
\newblock {\em Studia Logica}, 91(3):383--406.

\bibitem[Frankowski, 2004]{frankowski2004}
Frankowski, S. (2004).
\newblock Formalization of a plausible inference.
\newblock {\em Bulletin of the Section of Logic}, 33(1):41--52.

\bibitem[French and Ripley, 2018]{french2017valuations}
French, R. and Ripley, D. (2018).
\newblock Valuations: bi, tri, and tetra.
\newblock {\em Studia Logica}.
\newblock In press.

\bibitem[Gentzen, 1935]{gentzen1935investigations}
Gentzen, G. (1964 (1935)).
\newblock Investigations into logical deduction.
\newblock {\em American philosophical quarterly}, 1(4):288--306.
\newblock Original Publication in Mathematische Zeitschrift 39 (1): 176-210.

\bibitem[G\"odel, 1932]{godel1932intuitionistic}
G\"odel, K. (1932).
\newblock On the intuitionistic propositional calculus.
\newblock In Feferman, S., Dawson, J., Kleene, S., Moore, G., Solovay, R., and
  van Heijenoort, J., editors, {\em Kurt G\"odel, Collected Works, vol. 1},
  pages 223--225. Oxford University Press.

\bibitem[Humberstone, 1988]{humberstone1988heterogeneous}
Humberstone, L. (1988).
\newblock Heterogeneous logic.
\newblock {\em Erkenntnis}, 29(3):395--435.

\bibitem[Jansana, 2016]{sep-logic-algebraic-propositional}
Jansana, R. (2016).
\newblock Algebraic propositional logic.
\newblock In Zalta, E.~N., editor, {\em The Stanford Encyclopedia of
  Philosophy}. Metaphysics Research Lab, Stanford University, winter 2016
  edition.

\bibitem[Lahav and Zohar, 2018]{lahav2018subformula}
Lahav, O. and Zohar, Y. (2018).
\newblock From the subformula property to cut-admissibility in propositional
  sequent calculi.
\newblock {\em Journal of Logic and Computation}, 28(6):1341--1366.

\bibitem[Malinowski, 1990]{malinowski:q}
Malinowski, G. (1990).
\newblock Q-consequence operation.
\newblock {\em Reports on mathematical logic}, 24:49--54.

\bibitem[Marcos, 2009]{marcos2009non}
Marcos, J. (2009).
\newblock What is a non-truth-functional logic?
\newblock {\em Studia Logica}, 92(2):215.

\bibitem[Ripley, 2017]{ripley2017transitivity}
Ripley, D. (2017).
\newblock On the `transitivity' of consequence relations.
\newblock {\em Journal of Logic and Computation}, 28(2):433--450.

\bibitem[Scott, 1974]{scott1974completeness}
Scott, D. (1974).
\newblock Completeness and axiomatizability in many-valued logic.
\newblock In {\em Proceedings of the Tarski symposium}, volume~25, pages
  411--436. American Mathematical Society, Providence.

\bibitem[Shoesmith and Smiley, 1978]{shoesmith1978multiple}
Shoesmith, D.~J. and Smiley, T.~J. (1978).
\newblock {\em Multiple-conclusion logic}.
\newblock CUP Archive.

\bibitem[Shramko and Wansing, 2011]{shramko2011truth}
Shramko, Y. and Wansing, H. (2011).
\newblock {\em Truth and falsehood: An inquiry into generalized logical
  values}, volume~36.
\newblock Springer, Trends in Logic.

\bibitem[Smith, 2008]{smith:vdt}
Smith, N. J.~J. (2008).
\newblock {\em Vagueness and Degrees of Truth}.
\newblock Oxford University Press, Oxford.

\bibitem[Surma, 1982]{surma1982origin}
Surma, S.~J. (1982).
\newblock On the origin and subsequent applications of the concept of the
  {L}indenbaum algebra.
\newblock {\em Studies in Logic and the Foundations of Mathematics},
  104:719--734.

\bibitem[Suszko, 1975]{suszko1975remarks}
Suszko, R. (1975).
\newblock Remarks on{ \L}ukasiewicz's three-valued logic.
\newblock {\em Bulletin of the Section of Logic}, 3-4:87--89.

\bibitem[Suszko, 1977]{suszko1977fregean}
Suszko, R. (1977).
\newblock The {F}regean axiom and {P}olish mathematical logic in the 1920s.
\newblock {\em Studia Logica}, 36(4):377--380.

\bibitem[Tarski, 1930]{tarski1930}
Tarski, A. (1930).
\newblock On some fundamental concepts of metamathematics.
\newblock In Corcoran, J., editor, {\em Logic, Semantics, Metamathematics},
  pages 30--37.

\bibitem[Tsuji, 1998]{tsuji1998many}
Tsuji, M. (1998).
\newblock Many-valued logics and {S}uszko's thesis revisited.
\newblock {\em Studia Logica}, 60(2):299--309.

\bibitem[Wansing and Shramko, 2008]{wansing2008suszko}
Wansing, H. and Shramko, Y. (2008).
\newblock Suszko's thesis, inferential many-valuedness, and the notion of a
  logical system.
\newblock {\em Studia Logica}, 88(3):405--429.

\bibitem[W{\'o}jcicki, 1973]{wojcicki1973matrix}
W{\'o}jcicki, R. (1973).
\newblock Matrix approach in methodology of sentential calculi.
\newblock {\em Studia Logica}, 32(1):7--37.

\bibitem[W{\'o}jcicki, 1988]{wojcicki1988theory}
W{\'o}jcicki, R. (1988).
\newblock {\em Theory of logical calculi: basic theory of consequence
  operations}, volume 147.
\newblock Springer.

\end{thebibliography}

\end{document}